\documentclass[10pt]{article}
\usepackage{amsfonts}
\usepackage{mathrsfs}
\usepackage{amsmath,amsthm,amsfonts,amssymb,epsfig}

\newtheorem{thm}{Theorem}[section]

\newtheorem{lem}[thm]{Lemma}

\theoremstyle{definition}

\theoremstyle{remark}

\numberwithin{equation}{section}

\newcommand\be{\begin{equation}}
\newcommand\ee{\end{equation}}
\newcommand\bes{\begin{eqnarray}}
\newcommand\ees{\end{eqnarray}}
\newcommand\bess{\begin{eqnarray*}}
\newcommand\eess{\end{eqnarray*}}

\def\theequation{\arabic{section}.\arabic{equation}}
\begin{document}
\title{ Carleman estimates for a stochastic degenerate parabolic equation and applications to null controllability and an inverse random source problem}
\author{Bin Wu$^a$\footnote{Corresponding author. email: binwu@nuist.edu.cn}\quad  Qun Chen$^a$\quad Zewen Wang$^b$  \\
$ ^a$School of Mathematics and Statistics\\ Nanjing University of
Information Science
and Technology, \\
Nanjing 210044, P. R. China\\
$^b$School of Science\\
 East China University of Technology\\
  Nanchang 330013, P. R. China\\
}

\maketitle

\begin{abstract}
 In this paper, we establish two Carleman estimates for a stochastic degenerate parabolic equation. The first one is for the backward stochastic degenerate parabolic equation with singular weight function. Combining this Carleman estimate and an approximate argument, we prove the null controllability of the forward stochastic degenerate parabolic equation with the gradient term. The second one is for the forward stochastic degenerate parabolic equation with regular weighted function, based on which we obtain the Lipschitz stability for an inverse problem of determining a random source depending only on time in the forward stochastic degenerate parabolic equation.

 \vskip 0.3cm

{\bf AMS Subject Classifications:} 93B05, 93B07, 35R30

\vskip 0.3cm

{\bf Keywords:}  Stochstic degenerate parabolic equation, Carleman estimate, null controllability, inverse random source problem.

\end{abstract}
\section{Introduction}
\setcounter{equation}{0}

Let $(\Omega, \mathcal{F}, \{\mathcal{F}_t\}_{t\geq0}, \mathbb P)$ be a complete filtered probability space on which a one-dimensional standard Brownian motion $\{B(t)\}_{t\geq0}$ is defined such that $\{\mathcal F_t\}_{t\geq 0}$ is the natural filtration generated by $B(\cdot)$, augmented by all the $\mathbb P$-null sets in $\mathcal{F}$. Let $I=(0, 1)$ and $Q_T=I\times(0, T)$. Then, we consider the following forward stochastic  degenerate parabolic equation
\begin{align}\label{1.1}\left\{
\begin{array}{ll}
{\rm d}y-\left(x^\alpha y_x\right)_x{\rm d}t=f{\rm d}t+F{\rm d}B(t),&(x,t)\in Q_T,\\
y(1,t)=0,& t\in (0,T),\\
 {\rm and}\ \left\{\begin{array}{ll}
y(0,t)=0&{\rm for}\ \alpha\in (0,1),\\
\big(x^\alpha y_x\big)(0,t)=0 &{\rm for}\ \alpha\in [1,2),
\end{array} \right.
& t\in (0,T),\\
 y(x,0)=y_0(x),&x\in I,
 \end{array}\right.
\end{align}
and the following backward stochastic degenerate parabolic equation
\begin{align}\label{1.2}\left\{
\begin{array}{ll}
{\rm d}y+\left(x^\alpha y_x\right)_x{\rm d}t=f{\rm d}t+F{\rm d}B(t),&(x,t)\in Q_T,\\
y(1,t)=0,& t\in (0,T),\\
 {\rm and}\ \left\{\begin{array}{ll}
y(0,t)=0&{\rm for}\ \alpha\in (0,1),\\
\big(x^\alpha y_x\big)(0,t)=0 &{\rm for}\ \alpha\in [1,2),
\end{array} \right.
& t\in (0,T),\\
 y(x,T)=y_T(x),&x\in I.
 \end{array}\right.
\end{align}
Obviously, the equation is degenerate at the left-end point $x=0$.

The main objective of this paper is to obtain Carleman estimates for backward/forward stochastic degenerate
equations. As applications, we then apply these Carleman estimates to study a null controllability problem and an inverse random source problem. More precisely, for given subdomain $\omega=(x_1,x_2)$ such that $0<x_1<x_2<1$, we consider the following two problems:

\vspace{2mm}

{\noindent\bf Null controllability.} Find a pair control $(g,G)$ such that the solution $y$ of the following stochastic degenerate parabolic with the gradient term:
\begin{align}\label{1.3}\left\{
\begin{array}{ll}
{\rm d}y-\left(x^\alpha y_x\right)_x{\rm d}t=(a y_{x}+b y+g\mathbf 1_\omega){\rm d}t+(cy+G){\rm d}B(t),&(x,t)\in Q_T,\\
 y(0,t)=y(1,t)=0,& t\in (0,T),\\
 y(x,0)=y_0(x),&x\in I,
 \end{array}\right.
\end{align}
satisfies \begin{align*}y(x,T)=0, \quad x\in I,\ \mathbb P-a.s.,\end{align*}
where  $\alpha\in \big(0,\frac{1}{2}\big)$ and ${\mathbf 1}_{\omega}$ is the characteristic
function of the set $\omega$.

\vspace{2mm}

{\noindent\bf Remark 1.1.}\ For deterministic case,  [\ref{Wang2014}] pointed out that the restriction $\alpha\in \big(0,\frac{1}{2}\big)$ was optimal for establishing the Carleman
estimate under $a\in L^\infty(Q_T)$ according to the methods as [\ref{Alabau2006},\ref{Cannarsa2008}]. [\ref{Du2019}] gave a further explanation about this restriction. In other words,  $\alpha\in \big(0,\frac{1}{2}\big)$ is essentially caused by the tool used to prove the null controllability, i.e. Carleman estimates,  and is also  the best result based on the method used in this paper.

\vspace{2mm}

{\noindent\bf Inverse random source problem.} Determine $h(t)$ in the following stochastic degenerate parabolic equation:
\begin{align}\label{1.4}\left\{
\begin{array}{ll}
{\rm d}y-\left(x^\alpha y_x\right)_x{\rm d}t= h(t)r(x,t){\rm d}B(t),&(x,t)\in Q_T,\\
  y(1,t)=0,& t\in (0,T),\\
 {\rm and}\ \left\{\begin{array}{ll}
y(0,t)=0&{\rm for}\ \alpha\in (0,1),\\
\big(x^\alpha y_x\big)(0,t)=0 &{\rm for}\ \alpha\in [1,2),
\end{array} \right.
& t\in (0,T),\\
 y(x,0)=y_0(x),&x\in I,
 \end{array}\right.
\end{align}
by the observation data
\begin{align*}
 y|_{\omega_T}\quad {\rm and}\quad y(x,T).
\end{align*}

Carleman estimate is an important tool to study null controllability and inverse problems, which is a weighted estimate for a solution to a partial differential equation. There are rich references on Carleman estimates for deterministic partial differential equations, see [\ref{Bellassoued2015},\ref{Iman3},\ref{Isakov2},\ref{Klib3},\ref{Klib2},\ref{Rousseau},\ref{Wu2019},\ref{Wu2018},\ref{Yamamoto1}]. In recently years, Carleman estimates for stochastic partial differential equations are getting more and more attention. We refer to [\ref{Barbu2003},\ref{TangSIAM2009}] for stochastic
parabolic equation, [\ref{ZhangSIAM2008}] for stochastic hyperbolic equation, [\ref{Gao2014}] for stochastic Korteweg-de Vries equation, [\ref{Fu2017}] for stochastic complex Ginzburg-Landau equations and so on. To the best of our knowledge, there is only one paper about Carleman estimates for stochastic degenerate equations [\ref{LiuSIAM}], in which  the global Carleman estimates for some forward and backward stochastic degenerate parabolic equations were established and then were applied to an
insensitizing control problem.

 One  successful application of Carleman estimate in stochastic partial differential equations is to study related control problems for various mathematical models with stochastic effect [\ref{Gao2018},\ref{Liu2018MCRF},\ref{Lv2013},\ref{TangSIAM2009},\ref{Yan2018JMAA}].  As for null controllability for the deterministic degenerate equations, we refer to [\ref{Alabau2006},\ref{Cannarsa2008}] for degenerate parabolic equation, [\ref{Du2019},\ref{Wang2014}] for degenerate parabolic equation with the gradient terms, [\ref{Benhassi2013},\ref{Cannarsa2009-1},\ref{Du2018},\ref{Fadili2017}] for coupled degenerate systems and so on. On the other hand, there are few work on inverse problems for stochastic partial differential equations. We refer to [\ref{Lv2012}]  the uniqueness of an inverse source problem for the stochastic parabolic equation. An inverse source problem of
determining two kinds of sources simultaneously for a stochastic wave equation was studied in [\ref{Yan2015IP}]. Global uniqueness of an inverse problem of simultaneously determining random source and initial data for the stochastic hyperbolic equation in [\ref{Lv2015}]. This method then was extended to stochastic Euler-Bernoulli beam equation [\ref{Yuan2017JMAA}]. As for applications of regularization techniques in the numerical methods for inverse random source problems, we refer to [\ref{Bao2010}] or [\ref{Bao2013}].

In this paper, we first focus on Carleman estimates for stochastic degenerate parabolic equations. More precisely, we will prove two Carleman estimates for backward/forward stochastic degenerate parabolic equation, respectively with singular/regular weight functions.  We apply the first Carleman estimate with singular weight function to study the null controllability for stochastic degenerate parabolic equation with the gradient term (\ref{1.3}),  in whose proof we only assume that the coefficient  of the first
order term $a\in L^\infty_\mathcal F(0,T;L^\infty(I))$. Since the equation is degenerate, we could not apply directly the Carleman estimate to absorb the first order term, if $a\in L^\infty_\mathcal F(0,T;L^\infty(I))$.   To overcome this difficulty, we have to improve this Carleman estimate by using the method in [\ref{Imanuvilov2001},\ref{Du2019}] for deterministic differential equations, also see [\ref{Liu2014ECOCV}] for stochastic differential equations.  For this reason,
we only obtain the null controllability result for $\alpha\in
(0,\frac{1}{2})$. On the other hand, unlike the deterministic counterparts, the solution of a stochastic differential equation is not differentiable with respect to
time variable, Carleman estimate with singular weight function could not be applied
to inverse random source problem. Hence we would like to borrow some ideas from [\ref{Lv2015}] to prove the second Carleman estimate with regular weight function.  Applying this Carleman estimate, we obtain a Lipschitz stability for our inverse random source problem. In comparison with [\ref{LiuSIAM}], on one hand we release the power of $y$ on the left-hand side of Carleman estimate, which leads to that we can deal with a null controllability of stochastic degenerate equation with the first order term, see Theorem 4.2. On the other hand, since the weight function  in [\ref{LiuSIAM}] is singular in the Carleman estimate,  which could not be applied to study our inverse problem.

Throughout this paper, we denote by $L_\mathcal F^2(0,T)$ the space of all progressively measurable stochastic process $X$ such that $\mathbb E(\int_{0}^T|X|^2{\rm d}t)<\infty$.  For a Banach space $H$, we denote by  $L^2_\mathcal{F}(0, T; H)$ the Banach space consisting of all $H$-valued $\{\mathcal F_t\}_{t\geq0}$-adapted processes $X(\cdot)$ such that $\mathbb{E}(\|X(\cdot)\|_{L^2(0,T;H)}) $ $<\infty$, with the canonical norm; by $L^\infty_\mathcal{F}(0, T; H)$ the Banach space consisting of all $H$-valued $\{\mathcal F_t\}_{t\geq 0}$-adapted bounded processes; and by $L^2_\mathcal{F}(\Omega; C([0, T]; H))$ the Banach space consisting of all $H$-valued $\{\mathcal F_t\}_{t\geq0}$-adapted continuous processes $X(\cdot)$ such that $\mathbb{E}(\|X(\cdot)\|^2_{C([0,T];H)})<\infty$, with the canonical norm.

The remainder of this paper is organized as follows. In next section,
we prove the well-posedness of forward/backward stochastic degenerate parabolic equation with the first order term. In section 3, we show two Carleman estimates for backward/forward stochastic degenerate parabolic equations. In next two sections, based on these two Carleman estimates we study the null controllability and the inverse random source problem, respectively.

\section{Well-posedness}
\setcounter{equation}{0}

In this section, we use an approximate argument to prove the
well-posedness of the following stochastic degenerate parabolic equation
\begin{align}\label{2.1}\left\{
\begin{array}{ll}
{\rm d}y-\left(x^\alpha y_x\right)_x{\rm d}t=f{\rm d}t+F{\rm d}B(t),&(x,t)\in Q_T,\\
 y(1,t)=0,& t\in (0,T),\\
 {\rm and}\ \left\{\begin{array}{ll}
y(0,t)=0&{\rm for}\ \alpha\in (0,1),\\
\big(x^\alpha y_x\big)(0,t)=0 &{\rm for}\ \alpha\in [1,2),
\end{array} \right.
& t\in (0,T),\\
 y(x,0)=y_0(x),&x\in I,
 \end{array}\right.
\end{align}

To deal with degeneracy at $x=0$, we have to introduce following weighted space:
\begin{align*}
H_\alpha^1(I):=\left\{\zeta\in L^2(I)\ |\ x^{\frac{\alpha}{2}} \zeta\in L^2(I),\ \zeta(0)=\zeta(1)=0\right\}\quad {\rm for}\ \alpha\in (0,1)
\end{align*}
 and
\begin{align*}
H_\alpha^1(I):=\left\{\zeta\in L^2(I)\ |\ x^{\frac{\alpha}{2}} \zeta\in L^2(I),\ \zeta(1)=0\right\}\quad {\rm for}\ \alpha\in [1,2).
\end{align*}
We endow the space $H_\alpha^1(I)$ with the norm
\begin{align*}
\|\zeta\|^2_{H_\alpha^1}=\int_{I}\big(|\zeta|^2+x^{\alpha}|\zeta_x|^2\big){\rm d}x.
\end{align*}
Further we set
\begin{align*}
&\mathcal H^1= L^2_\mathcal F(\Omega;C([0,T];L^2(I)))\cap L^2_\mathcal F(0,T;H^1(I))\nonumber\\
&\mathcal H_\alpha^1=L^2_\mathcal F(\Omega;C([0,T];L^2(I)))\cap L^2_\mathcal F(0,T;H_\alpha^1(I)).
\end{align*}

{\noindent\bf Definition.}\ A stochastic process $y$ is said to be a weak solution of the forward stochastic degenerate parabolic equation (\ref{2.1}) if $y\in \mathcal H^1_\alpha$ and $y(0)=y_0$ in $I$, $\mathbb P-a.s.$ and it holds for all $\phi\in C_0^\infty(\overline I)$ that
\begin{align*}
&\int_{I}y(x,t)\phi(x){\rm d}x-\int_{I}y_0(x)\phi(x){\rm d}x+\int_{Q_t}x^{\alpha}y_{x}\phi_x{\rm d}x{\rm d}t\nonumber\\
=&\int_{Q_t}f\phi{\rm d}x{\rm d}t+\int_{Q_t}F\phi{\rm d}x{\rm d}B(t),\quad \mathbb P-a.s.
\end{align*}

\begin{thm} Let $f, F\in L^2_\mathcal F(0,T;L^2(I))$ and $y_0\in L^2(\Omega,\mathcal F_0,\mathbb P;L^2(I))$. Then (\ref{2.1}) admits a unique weak solution $y\in \mathcal H^1_\alpha$.
\end{thm}

{\noindent\bf Proof.}\ Letting $\varepsilon\in (0,1)$,
we consider the following nondegenerate approximate problem:
\begin{align}\label{2.2}\left\{
\begin{array}{ll}
{\rm d}y^\varepsilon-\left((x+\varepsilon)^\alpha y^\varepsilon_x\right)_x{\rm d}t=f{\rm d}t+F{\rm d}B(t),&(x,t)\in Q_T,\\
 y^\varepsilon(1,t)=0,& t\in (0,T),\\
 {\rm and}\ \left\{\begin{array}{ll}
y^\varepsilon(0,t)=0&{\rm for}\ \alpha\in (0,1),\\
\big((x+\varepsilon)^\alpha y^\varepsilon_x\big)(0,t)=0 &{\rm for}\ \alpha\in [1,2),
\end{array} \right.
& t\in (0,T),\\
 y^\varepsilon(x,0)=y^\varepsilon_0(x),&x\in I,
 \end{array}\right.
\end{align}
where
\begin{align}\label{1-2.3}
y^\varepsilon_0\rightarrow y_0\quad {\rm in}\ L^2(\Omega,\mathcal F_0,\mathbb P;L^2(I)).
\end{align}
Then by [\ref{Krylov1994}] or [\ref{Hu1991}], it is easy to check that (\ref{2.2}) admits a unique weak solution $y^\varepsilon\in \mathcal H^1$.

Now we prove a uniform estimate in $\varepsilon$ for $y^\varepsilon$:
\begin{align}\label{2.3}
&\mathbb E\sup_{t\in [0,T]}\int_{I}|y^\varepsilon|^2(x,t){\rm d}x+\mathbb E\int_{Q_T} (x+\varepsilon)^{\alpha}|y^\varepsilon_{x}|^2{\rm d}x{\rm d}t
\nonumber\\
\leq &C\mathbb E\int_{I}|y^\varepsilon_0|^2{\rm d}x+C\mathbb E\int_{Q_T} (|f|^2+|F|^2){\rm d}x{\rm d}t.
\end{align}
where $C$ is depending on $I,T$ and $\alpha$, but independent of $\varepsilon$.
By It\^{o} formula and the equation of $y^\varepsilon$, we obtain
\begin{align}\label{2.4}
&\int_{I}|y^\varepsilon|^2(x,t){\rm d}x+2\int_{Q_t} (x+\varepsilon)^{\alpha}|y^\varepsilon_{x}|^2{\rm d}x{\rm d}t
\nonumber\\
= &\int_{I}|y^\varepsilon_0|^2{\rm d}x+2\int_{Q_t} fy^\varepsilon{\rm d}x{\rm d}t+\int_{Q_t}|F|^2{\rm d}x{\rm d}t+2\int_{Q_t}Fy^\varepsilon{\rm d}x{\rm d}B(t).
\end{align}
By the Burkholder-Davis-Gundy inequality, we can obtain for any $\epsilon>0$ that
\begin{align}\label{1-2.6}
\mathbb E\left(\sup_{t\in[0,T]}\left|\int_{Q_t}Fy^\varepsilon{\rm d}x{\rm d}B(t)\right|\right)\leq \epsilon\mathbb E\left(\sup_{t\in [0,T]}\int_{I}|y^\varepsilon|^2(x,t){\rm d}x\right)+C(\epsilon)\mathbb E\int_{Q_T}|F|^2{\rm d}x{\rm d}t.
\end{align}
By (\ref{2.4}) and (\ref{1-2.6}), we have
\begin{align}\label{2.5}
&\mathbb E\sup_{t\in [0,T]}\int_{I}|y^\varepsilon|^2(x,t){\rm d}x\leq C\mathbb E\int_{I}|y_0|^2{\rm d}x+C\mathbb E\int_{Q_T}\big(|f|^2+|F|^2\big){\rm d}x{\rm d}t.
\end{align}
Moreover, it follows from the Burkholder-Davis-Gundy inequality and (\ref{2.4}) that
\begin{align}\label{2.6}
&\mathbb E\sup_{t\in [0,T]}\int_{I}|y^\varepsilon|^2(x,t){\rm d}x+\mathbb E\int_{Q_T} (x+\varepsilon)^{\alpha}|y^\varepsilon_{x}|^2{\rm d}x{\rm d}t
\nonumber\\
\leq &\mathbb E\int_{I}|y_0|^2{\rm d}x+C\mathbb E\int_{Q_T} |y^\varepsilon|^2{\rm d}x{\rm d}t+C\mathbb E\int_{Q_T} (|f|^2+|F|^2){\rm d}x{\rm d}t.
\end{align}
Substituting (\ref{2.5}) into (\ref{2.6}), we obtain (\ref{2.3}).

Similarly, we have for any $\varepsilon_1,\varepsilon_2\in (0,1)$ that
\begin{align}
&\mathbb E\sup_{t\in [0,T]}\int_{I}|y^{\varepsilon_1}-y^{\varepsilon_2}|^2(x,t){\rm d}x+\mathbb E\int_{Q_T} (x+\varepsilon)^{\alpha}|y^{\varepsilon_1}_{x}-y^{\varepsilon_2}_{x}|^2{\rm d}x{\rm d}t
\nonumber\\
\leq &C\mathbb E\int_{I}|y^{\varepsilon_1}_0-y^{\varepsilon_2}_0|^2{\rm d}x.
\end{align}
Therefore, $\{y^\varepsilon\}$ is a Cauchy sequence in $\mathcal H^1_\alpha$. Letting $\varepsilon\rightarrow0$, we find that (\ref{2.1}) admits  a weak solution $y\in \mathcal H_\alpha$ (the limit of $y^\varepsilon$ in $\mathcal H_\alpha$). The uniqueness of solution could be directly deduced from (\ref{2.3}). \hfill$\Box$

\vspace{2mm}

Next we consider the stochastic degenerate parabolic equation with gradient term:
\begin{align}\label{2.8}\left\{
\begin{array}{ll}
{\rm d}y-\left(x^\alpha y_x\right)_x{\rm d}t=(a y_{x}+b y+f){\rm d}t+(cy+F){\rm d}B(t),&(x,t)\in Q_T,\\
 y(0,t)=y(1,t)=0,& t\in (0,T),\\
 y(x,0)=y_0(x),&x\in I.
 \end{array}\right.
\end{align}
In comparison with (\ref{2.1}), the main difficulty is how to deal with the gradient term under $a\in L^\infty_\mathcal F(0,T;L^\infty(I))$. Due to degeneracy, we could not control this term directly. We apply the method in [\ref{Wang2014}] to overcome this difficulty.  Based on this reason, we only prove the well-posedness of (\ref{2.8}) when $\alpha\in (0,1)$.

\begin{thm} Let $\alpha\in (0,1)$, $a,b,c\in L^\infty_\mathcal F(0,T;L^\infty(I))$, $f,F\in L^2_\mathcal F(0,T;L^2(I))$ and $y_0\in L^2(\Omega,\mathcal F_0,\mathbb P;L^2(I))$. Then (\ref{2.8}) admits a unique weak solution $y\in \mathcal H^1_\alpha$.
\end{thm}

{\noindent\bf Proof.}\ We also use an approximate argument to prove this result. Let $y^\varepsilon\in \mathcal H^1$ be the unique solution of the following problem:\begin{align}\label{2.10}\left\{
\begin{array}{ll}
{\rm d}y^\varepsilon-\left((x+\varepsilon)^\alpha y^\varepsilon_x\right)_x{\rm d}t=(a y^\varepsilon_{x}+b y^\varepsilon+f){\rm d}t+(cy^\varepsilon+F){\rm d}B(t),&(x,t)\in Q_T,\\
 y^\varepsilon(0,t)=y^\varepsilon(1,t)=0,& t\in (0,T),\\
 y^\varepsilon(x,0)=y^\varepsilon_0(x),&x\in I,
 \end{array}\right.
\end{align}
where the sequence $\{y_0^\varepsilon\}$ satisfies (\ref{1-2.3}). Then similar to (\ref{2.4}), we have
\begin{align}\label{2.11}
&\int_{I}|y^\varepsilon|^2(x,t){\rm d}x+2\int_{Q_t} (x+\varepsilon)^{\alpha}|y^\varepsilon_{x}|^2{\rm d}x{\rm d}t
\nonumber\\
= &\int_{I}|y^\varepsilon_0|^2{\rm d}x+2\int_{Q_T}a y^\varepsilon_{x}y^\varepsilon{\rm d}x{\rm d}t+2\int_{Q_t} (b y^\varepsilon+f)y^\varepsilon{\rm d}x{\rm d}t+\int_{Q_t}|cy^\varepsilon+F|^2{\rm d}x{\rm d}t\nonumber\\
&+2\int_{Q_t}Fy^\varepsilon{\rm d}x{\rm d}B(t)\nonumber\\
\leq &\int_{I}|y^\varepsilon_0|^2{\rm d}x+2\int_{Q_T}a y^\varepsilon_{x}y^\varepsilon{\rm d}x{\rm d}t+\int_{Q_T} (|f|^2+|F|^2){\rm d}x{\rm d}t+C\int_{Q_t}|y^\varepsilon|^2{\rm d}x{\rm d}t\nonumber\\
&+2\int_{Q_t}Fy^\varepsilon{\rm d}x{\rm d}B(t).
\end{align}
By Young's inequality, we have
\begin{align}\label{2.12}
\int_{Q_t} ay^\varepsilon_xy^\varepsilon {\rm d}x{\rm d}t\leq& \frac{1}{4} \int_{Q_t} (x+\varepsilon)^\alpha |y^\varepsilon_x|^2{\rm d}x{\rm d}t+C\int_{Q_t}(x+\varepsilon)^{-\alpha}|y^\varepsilon|^2{\rm d}x{\rm d}t.
\end{align}
For a sufficiently small $\kappa>0$ we have
\begin{align}\label{2.13}
&\int_{Q_t}(x+\varepsilon)^{-\alpha}|y^\varepsilon|^2{\rm d}x{\rm d}t\nonumber\\
=&\int_{0}^{t}\int_0^\kappa(x+\varepsilon)^{-\alpha}|y^\varepsilon|^2{\rm d}x{\rm d}t+\int_{0}^{t}\int_\kappa^1(x+\varepsilon)^{-\alpha}|y^\varepsilon|^2{\rm d}x{\rm d}t\nonumber\\
\leq&\int_{0}^{t}\int_0^\kappa(x+\varepsilon)^{-\alpha}\left|\int_0^x y^\varepsilon_x(\zeta,t){\rm d}\zeta\right|^2{\rm d}x{\rm d}t+ \int_{0}^{t}\int_\kappa^1 (x+\varepsilon)^{-\alpha}|y^\varepsilon|^2{\rm d}x{\rm d}t\nonumber\\
\leq &\frac{1}{(1-\alpha)}\int_{0}^{t}\int_0^\kappa (x+\varepsilon)^{1-2\alpha}\left(\int_0^x(\zeta+\varepsilon)^\alpha |y^\varepsilon_x(\zeta,t)|^2{\rm d}\zeta\right){\rm d}x{\rm d}t+\kappa^{-\alpha}\int_{0}^{t}\int_I|y^\varepsilon|^2{\rm d}x{\rm d}t\nonumber\\
\leq & C_{\kappa,\alpha}\int_{Q_t} (x+\varepsilon)^\alpha |y^\varepsilon_x|^2{\rm d}x{\rm d}t+\kappa^{-\alpha}\int_{Q_t}|y^\varepsilon|^2{\rm d}x{\rm d}t
\end{align}
with $$C_{\kappa,\alpha}=\frac{(\kappa+\varepsilon)^{2-2\alpha}-\varepsilon^{2-2\alpha}}{(1-\alpha)(2-2\alpha)}.$$ For any $\varepsilon\in (0,1)$, we choose $\kappa$ sufficiently small such that $C_{\kappa,\alpha}<\frac{1}{4C}$. Then, from (\ref{2.12}) and (\ref{2.13}) we deduce that
\begin{align}\label{2.14}
\int_{Q_t} a y^\varepsilon_x y^\varepsilon{\rm d}x{\rm d}t\leq& \frac{1}{2}  \int_{Q_t} (x+\varepsilon)^\alpha |y^\varepsilon_x|^2{\rm d}x{\rm d}t+C\int_{Q_t}|y^\varepsilon|^2{\rm d}x{\rm d}t.
\end{align}

Substituting (\ref{2.14}) into (\ref{2.11}) and taking mathematical expectation yields that
\begin{align}
&\mathbb E\int_{I}|y^\varepsilon|^2(x,t){\rm d}x+\mathbb E\int_{Q_t} (x+\varepsilon)^{\alpha}|y^\varepsilon_{x}|^2{\rm d}x{\rm d}t
\nonumber\\
\leq &\mathbb E\int_{I}|y_0|^2{\rm d}x+C\mathbb E\int_{Q_T} (|f|^2+|F|^2){\rm d}x{\rm d}t+C\mathbb E\int_{Q_t} |y^\varepsilon|^2{\rm d}x{\rm d}t,
\end{align}
which implies
\begin{align}
&\mathbb E\sup_{t\in [0,T]}\int_{I}|y^\varepsilon|^2(x,t){\rm d}x\leq C\mathbb E\int_{I}|y_0|^2{\rm d}x+C\mathbb E\int_{Q_T}\big(|f|^2+|F|^2\big){\rm d}x{\rm d}t
\end{align}
by Gronwall inequality. The remainder of the proof is almost the same as the one in Theorem 2.1.
\hfill$\Box$

\vspace{2mm}

{\noindent\bf Remark 2.1.}\ If $a$ has the decomposition $a=x^{\frac{\alpha}{2}}\tilde a$ with some $\tilde a\in L^\infty_\mathcal F(0,T;L^\infty(I))$ as in [\ref{Du2018}], the term of $a y^\varepsilon_x y^\varepsilon$ could be absorbed directly by the terms on the left-hand side of (\ref{2.11}). Then we can obtain the well-posedness  of (\ref{2.8}) for all $\alpha\in (0,2)$. Or if $ a\in L^\infty_\mathcal F(0,T;W^{1,\infty}(I))$, we also obtain the well-posedness for $\alpha\in (0,2)$.

\section{Carleman estimates for stochastic degenerate parabolic equation}
\setcounter{equation}{0}

In this section, we will show two Carleman estimates for stochastic degenerate parabolic equations. One is for the backward stochastic degenerate parabolic equation.
We will apply this Carleman estimate to prove the null controllability result for the forward stochastic degenerate parabolic equation with the gradient term. So that we use a singular weight function in this Carleman estimate. The other one is for the
forward stochastic degenerate parabolic equation, which will be used to study our inverse random source problem.  Unlike the deterministic case, we could not differentiate the stochastic equation with respect to time. For this reason, in order to prove the Lipschitz stability of our inverse problem, we have to introduce a regular weight function to put the term of unknown random source on the left-hand side of Carlmen estimate.


\subsection{Carleman estimate for backward stochastic degenerate equation}

We first introduce some weight functions. For
$\omega=(x_1,x_2)$, we choose $\omega^{(i)}:=(x_1^{(i)},x_2^{(i)})
(i=1,2)$ such that
$\omega^{(2)}\Subset\omega^{(1)}\Subset\omega$.
Let  $\chi\in C^2(\overline I)$ be a cut-off function such that
$0\leq \chi(x)\leq 1$ for $x\in  I$, $\chi(x)\equiv1$ for $x\in
(0,x_1^{(2)})$ and $\chi(x)\equiv0$ for $x\in (x_2^{(2)},1)$.  For a suitable positive constant $\beta$, we introduce
\begin{align*} \eta_1(x)=(x+\varepsilon)^{\beta}-\varepsilon^\beta,\quad x\in I,
\end{align*} and
$\eta_2\in C^2(\overline I)$ such that
\begin{align*}
\eta_2(x)>0,\quad x\in I,\quad \eta_2(0)=\eta_2(1)=0\quad{\rm
and}\quad |\eta_{2,x}(x)|>0,\quad  x\in
\overline{I\setminus\omega^{(1)}}
\end{align*}
and
\begin{align*}
 \eta_1(x)=\eta_2(x),\quad x\in (x^{(2)}_1,x_2^{(2)}).
\end{align*}
 Let us define
\begin{align*}
&\xi(t)=\frac{1}{t^2(T-t)^2},\quad
\psi_i(x)=e^{\lambda\eta_i(x)}-e^{2\lambda M},\quad \phi_i(x)=e^{\lambda\eta_i(x)},\quad i=1,2.
\end{align*}
where $\lambda$ is a positive parameter and $M$ is a sufficiently large constant such that $$M\geq \max\big\{\|\eta_1\|_{C(\overline\Omega)}, \|\eta_2\|_{C(\overline\Omega)}\big\}.$$
Now we introduce weight function in the first Carleman estimate
\begin{align*}
&\varphi(x,t)=\chi(x)\varphi_1(x,t)+(1-\chi(x))\varphi_2(x,t),\quad (x,t)\in Q_T
\end{align*}
with
 \begin{align*}\varphi_i(x,t)=\psi_i(x)\xi(t),\quad i=1,2.\end{align*}

We easily see that
\begin{align}\label{3.1}
\varphi(x,t)=\left\{\begin{array}{ll}\varphi_1(x,t),&(x,t)\in (0,x_1^{(2)})\times (0,T),\\
\varphi_1(x,t)=\varphi_2(x,t),&(x,t)\in (x_1^{(2)},x_2^{(2)})\times (0,T),\\
\varphi_2(x,t),&(x,t)\in (x_2^{(2)},1)\times (0,T).\\
\end{array}\right.
\end{align}

In order to deal with degeneracy, we first prove the following uniform Carleman estimate in $\varepsilon$.

\begin{thm} Let $\alpha\in \left(0,2\right)$, $f_1\in L^2_\mathcal F(0,T;L^2(I))$, $F_1\in L^2_\mathcal F(0,T; H^1(I))$, $u_T^\varepsilon\in L^2(\Omega,\mathcal F_T,$ $\mathbb P; L^2(I))$ and $\beta$ such that
\begin{align}\label{3.2}\left\{\begin{array}{ll}
1<\beta\leq 2-\alpha,& \alpha\in (0,1),\\
\beta=1,&\alpha=1,\\
\beta_0<\beta\leq 2-\alpha,&\alpha\in (1,2),
\end{array}
\right.
\end{align}
with $$\beta_0=\max\left\{0,3-2\alpha,1-\frac{\alpha}{2},\frac{14-9\alpha+\sqrt{17\alpha^2-44\alpha+36}}{8}\right\}.$$ Then for any $\varepsilon\in (0,1)$, there exist positive constants $\lambda_1=\lambda_1(\omega,I,T,\alpha,M)$, $s_1=s_1(\omega,I,T$, $\alpha,M,\lambda)$ and $C_1=C_1(\omega,I,T,\alpha,M)$, $C_2=C_2(\omega,I,T,\alpha,$ $M$, $\lambda)$
such that
\begin{align}\label{3.3}
&\mathbb E\int_{Q_T}s^3\xi^3(x+\varepsilon)^{2\alpha+3\beta-4}|u^\varepsilon|^2e^{2s\varphi}{\rm d}x{\rm
d}t+\mathbb E\int_{Q_T}s\xi(x+\varepsilon)^{2\alpha+\beta-2}
|u^\varepsilon_x|^2e^{2s\varphi}{\rm d}x{\rm d}t\nonumber\\
\leq& C_1\mathbb E\int_{Q_T} |f_1|^2e^{2s\varphi}{\rm d}x{\rm
d}t+C_2(\lambda)\mathbb E\int_{Q_T}
s^{2}\xi^{2}|F_1|^2e^{2s\varphi}{\rm d}x{\rm
d}t\nonumber\\
&+C_2(\lambda)\mathbb E\int_{\omega_T} s^3\xi^3|u^\varepsilon|^2 e^{2s\varphi}{\rm
d}x{\rm d}t+C_2(\lambda)\mathbb E\int_0^T s^2\left[(x+\varepsilon)^{\alpha+\beta-1}\xi^3 |u^\varepsilon|^2e^{2s\varphi}\right]_{x=0}{\rm d}t
\end{align}
for all $\lambda\geq \lambda_1$, $s\geq s_1$ and all $u^\varepsilon\in \mathcal H^1$ satisfying
\begin{align}\label{3.4}\left\{
\begin{array}{ll}
{\rm d}u^\varepsilon+\left((x+\varepsilon)^\alpha u^\varepsilon_x\right)_x{\rm d}t=f_1{\rm d}t+F_1{\rm d}B(t),&(x,t)\in Q_T,\\
u^\varepsilon(1,t)=0,& t\in (0,T),\\
 {\rm and}\ \left\{\begin{array}{ll}
 u^\varepsilon(0,t)=0 & {\rm for}\ \alpha\in (0,1),\\
 ((x+\varepsilon)^\alpha u^\varepsilon_x)(0,t)=0 & {\rm for}\ \alpha\in [1,2),
 \end{array}\right. &t\in (0,T),\\
 u^\varepsilon(x,T)=u^\varepsilon_T(x),&x\in I.
 \end{array}\right.
\end{align}
\end{thm}

%

{\noindent\bf Remark 3.1.}\ Given any $\varepsilon\in (0,1)$, the equation (\ref{3.4}) is not degenerate. Therefore, the regularity $u\in \mathcal H^1$ we assumed in Theorem 3.1 is reasonable.

\vspace{2mm}

{\noindent\bf Remark 3.2.}\ For $\alpha\in (1,2)$, $\beta\in (\beta_0,2-\alpha)$ is nonempty.

\vspace{2mm}

Letting $\varepsilon\rightarrow 0$ in Theorem 3.1, we can obtain the following Carleman estimate:

\begin{thm} Let $\alpha\in \left(0,2\right)$, $f_1\in L^2_\mathcal F(0,T;L^2(I))$, $F_1\in L^2_\mathcal F(0,T; H^1(I))$, $u_T\in L^2(\Omega,\mathcal F_T,$ $\mathbb P; L^2(I))$. Then for any $\varepsilon\in (0,1)$, there exist positive constants $\lambda_1=\lambda_1(\omega,I,T,\alpha,M)$, $s_1=s_1(\omega,I,T$, $\alpha,M,\lambda)$ and $C_1=C_1(\omega,I,T,\alpha,M)$, $C_2=C_2(\omega,I,T,\alpha,$ $M$, $\lambda)$
such that
\begin{align}\label{1-3.5}
&\mathbb E\int_{Q_T}s^3\xi^3x^{2-\alpha}|u|^2e^{2s\varphi}{\rm d}x{\rm
d}t+\mathbb E\int_{Q_T}s\xi x^{\alpha}
|u_x|^2e^{2s\varphi}{\rm d}x{\rm d}t\nonumber\\
\leq& C_1\mathbb E\int_{Q_T} |f_1|^2e^{2s\varphi}{\rm d}x{\rm
d}t+C_2(\lambda)\mathbb E\int_{Q_T}
s^{2}\xi^{2}|F_1|^2e^{2s\varphi}{\rm d}x{\rm
d}t\nonumber\\
&+C_2(\lambda)\mathbb E\int_{\omega_T} s^3\xi^3|u^\varepsilon|^2 e^{2s\varphi}{\rm
d}x{\rm d}t
\end{align}
for all $\lambda\geq \lambda_1$, $s\geq s_1$ and all $u\in \mathcal H^1_\alpha$ satisfying
\begin{align}\label{1-3.6}\left\{
\begin{array}{ll}
{\rm d}u+\left(x^\alpha u_x\right)_x{\rm d}t=f_1{\rm d}t+F_1{\rm d}B(t),&(x,t)\in Q_T,\\
u(1,t)=0,& t\in (0,T),\\
 {\rm and}\ \left\{\begin{array}{ll}
 u(0,t)=0 & {\rm for}\ \alpha\in (0,1),\\
 (x^\alpha u_x)(0,t)=0 & {\rm for}\ \alpha\in [1,2),
 \end{array}\right. &t\in (0,T),\\
 u(x,T)=u_T(x),&x\in I.
 \end{array}\right.
\end{align}
\end{thm}

{\noindent\bf Proof.}\ Letting $u^\varepsilon_T\rightarrow u_T$ in $L^2(\Omega,\mathcal F_T,\mathbb P;L^2(I))$, we easily see that $u^\varepsilon\rightarrow u$ in $\mathcal H^1_\alpha$. Then for a.e. $t\in (0,T)$, we have $u(\cdot,t)\in H^1_\alpha$, which implies that $xu^2\in W^{1,1}(I)$ and $x u^2\rightarrow 0$ as $x\rightarrow 0$ by Lemma 3.5 in [\ref{Cannarsa2008}]. Then we set $\beta=2-\alpha$  and choose $\varepsilon\rightarrow 0$ in (\ref{3.3}) to obtain (\ref{1-3.5}). \hfill$\Box$

\vspace{2mm}

 To prove Theorem 3.1, we need the following two lemmas. One is Hardy-Poincar\'{e} inequality [\ref{Opic1990}]. The detailed proof could be found in [\ref{Cannarsa2005}] or [\ref{Alabau2006}]. The other one is the Cacciopoli inequality for the stochastic parabolic equation, whose proof is detailed in Appendix and omitted here.

\begin{lem} Let $\gamma\in [0,1)\cup (1,2]$ and $z\in H_{0}^1(G)$. Then for any
$\varepsilon\in (0,1)$, we have
\begin{align}\label{3.5}
\int_I(x+\varepsilon)^{-\gamma}|z|^2{\rm d}x\leq
\frac{4}{(\gamma-1)^2}\int_I
(x+\varepsilon)^{2-\gamma}|z_x|^2{\rm d}x.
\end{align}
\end{lem}

\begin{lem} Let  $f_1\in L^2_\mathcal F(0,T;L^2(I))$ and $F_1\in L^2_\mathcal F(0,T; L^2(I))$. Then there exist positive constants $C_3=C_3(\omega,I,T,\alpha,M)$ and $C_4=C_4(\omega,I,T,\alpha,M,\lambda)$ such that the solution $u\in \mathcal H$ of  the backward stochastic degenerate parabolic equation (\ref{3.4})
 satisfies
 \begin{align}\label{3.7}
\mathbb E \int_{\omega^{(1)}_T} \xi|u^\varepsilon_x|^2e^{2s\varphi}{\rm d}x{\rm d}t\leq& C_4(\lambda)\mathbb E \int_{\omega_T}s^2\xi^3|u^\varepsilon|^2e^{2s\varphi}{\rm d}x{\rm d}t+C_3\mathbb E\int_{Q_T} s^{-2}|f_1|^2e^{2s\varphi}{\rm d}x{\rm d}t\nonumber\\
&+C_3\mathbb E \int_{Q_T} \xi|F_1|^2e^{2s\varphi}{\rm d}x{\rm d}t.
\end{align}
\end{lem}

Now we prove Theorem 3.1.

\vspace{2mm}

\noindent{\bf Proof of Theorem 3.1.}\   We split the proof into the following four steps.

{\em Step 1. A weighted identity for  backward stochastic degenerate parabolic operator.}

Let $l_1=s\varphi_1$, $\theta_1=e^{l_1}$ and $U=\theta_1 u^\varepsilon$. Then we have
\begin{align}\label{3.8}
\theta_1 \left[{\rm d}u^\varepsilon+\big((x+\varepsilon)^\alpha u^\varepsilon_x\big)_x{\rm d}t\right]=I_1+I_2{\rm d}t
\end{align}
with
\begin{align}\label{1-3.8}\left\{\begin{array}{ll}
U(1,t)=0,& t\in (0,T),\\
{\rm and}\ \left\{\begin{array}{ll}
U(0,t)=0& {\rm for}\ \alpha\in (0,1),\\
((x+\varepsilon)^\alpha U_x)(0,t)=\big((x+\varepsilon)^{\alpha}l_{1,x}U\big)(0,t)& {\rm for}\ \alpha\in [1,2),
\end{array}
\right. &t\in (0,T)\\
U(x,0)=U(x,T)=0,& x\in \overline I,\end{array}\right.
\end{align}
where
\begin{align*}
&I_1={\rm d}U-2(x+\varepsilon)^\alpha l_{1,x}U_x{\rm d}t-\big((x+\varepsilon)^\alpha l_{1,x}\big)_x U{\rm d}t,\\
&I_2= \big((x+\varepsilon)^\alpha U_{x}\big)_x+(x+\varepsilon)^{\alpha}l_{1,x}^2 U-l_{1,t} U.
\end{align*}
Hence,
\begin{align}\label{3.9}
\theta_1 I_2\left[{\rm d}u^\varepsilon+\big((x+\varepsilon)^\alpha u^\varepsilon_x\big)_x{\rm d}t\right]=I_1I_2+|I_2|^2{\rm d}t.
\end{align}
Now we deal with the term $I_1I_2$. By applying It\^{o} formula, we obtain
\begin{align}\label{3.10}
I_2{\rm d}U
=&\big[(x+\varepsilon)^\alpha U_x{\rm d}U\big]_x-\frac{1}{2}{\rm d}\big[(x+\varepsilon)^\alpha U_x^2\big]+\frac{1}{2}(x+\varepsilon)^\alpha\big({\rm d}U_x\big)^2\nonumber\\
&+\frac{1}{2}{\rm d}\big[(x+\varepsilon)^\alpha l_{1,x}^2 U^2\big]-(x+\varepsilon)^\alpha l_{1,x}l_{1,xt}U^2{\rm d}t-\frac{1}{2}(x+\varepsilon)^\alpha l_{1,x}^2({\rm d}U)^2\nonumber\\
&-\frac{1}{2}{\rm d} (l_{1,t}U^2)+\frac{1}{2}l_{1,tt} U^2{\rm d}t+\frac{1}{2}l_{1,t}({\rm d}U)^2.
\end{align}
On the other hand, a direct calculation yields
\begin{align}\label{3.11}
&I_2\Big[-2(x+\varepsilon)^\alpha l_{1,x}U_x{\rm d}t-\big((x+\varepsilon)^\alpha l_{1,x}\big)_x U{\rm d}t\Big]\nonumber\\
=&-\big[(x+\varepsilon)^{2\alpha}l_{1,x} U_x^2\big]_x{\rm d}t+(x+\varepsilon)^{2\alpha}l_{1,xx} U_x^2{\rm d}t\nonumber\\
&-\big[(x+\varepsilon)^{2\alpha}l_{1,x}^3 U^2\big]_x{\rm d}t+\big[2\alpha(x+\varepsilon)^{2\alpha-1}l_{1,x}^3+3(x+\varepsilon)^{2\alpha}l_{1,x}^2 l_{1,xx}\big]U^2{\rm d}t\nonumber\\
&+\big[(x+\varepsilon)^\alpha l_{1,x}l_{1,t} U^2\big]_x{\rm d}t-\big[(x+\varepsilon)^\alpha l_{1,x}l_{1,t}\big]_x U^2{\rm d}t\nonumber\\
&-\Big[\big((x+\varepsilon)^\alpha l_{1,x}\big)_x(x+\varepsilon)^\alpha U U_x\Big]_x{\rm d}t+\big[(x+\varepsilon)^\alpha l_{1,x}\big]_{xx} (x+\varepsilon)^\alpha UU_x{\rm d}t\nonumber\\
&+\big[(x+\varepsilon)^\alpha l_{1,x}\big]_{x}(x+\varepsilon)^\alpha U_x^2{\rm d}t-\big[(x+\varepsilon)^{2\alpha} l_{1,x}^2 l_{1,xx}+\alpha (x+\varepsilon)^{2\alpha-1} l_{1,x}^3\big]U^2{\rm d}t\nonumber\\
&+\big[(x+\varepsilon)^\alpha l_{1,x}\big]_x l_{1,t} U^2{\rm d}t.
\end{align}

Therefore, by (\ref{3.9})-(\ref{3.11}), we obtain the following weighted identity
\begin{align}\label{1-3.11}
&\theta_1 I_2\big[{\rm d}u^\varepsilon+\big((x+\varepsilon)^\alpha u^\varepsilon_x\big)_x{\rm d}t\big]=|I_2|^2{\rm d}t+K_1{\rm d}t+(K_2)_x+{\rm d}K_3+K_4,
\end{align}
where
\begin{align*}
K_1=&\left[\alpha(x+\varepsilon)^{2\alpha-1}l_{1,x}^3+2(x+\varepsilon)^{2\alpha}l_{1,x}^2 l_{1,xx}\right]U^2\\
&+\big[2(x+\varepsilon)^{2\alpha}l_{1,xx}+\alpha(x+\varepsilon)^{2\alpha-1}l_{1,x}\big]U_x^2+(x+\varepsilon)^\alpha\big[(x+\varepsilon)^\alpha l_{1,x}\big]_{xx}  UU_x\nonumber\\
&-(x+\varepsilon)^\alpha l_{1,x}l_{1,xt}U^2+\frac{1}{2}l_{1,tt} U^2-\big[(x+\varepsilon)^\alpha l_{1,x}l_{1,t}\big]_x U^2\nonumber\\
&+\big[(x+\varepsilon)^\alpha l_{1,x}\big]_x l_{1,t} U^2,
\end{align*}
\begin{align*}
K_2=&(x+\varepsilon)^\alpha U_x{\rm d}U-(x+\varepsilon)^{2\alpha} l_{1,x} U_x^2{\rm d}t-(x+\varepsilon)^{2\alpha}l_{1,x}^3 U^2{\rm d}t\\
&+(x+\varepsilon)^\alpha l_{1,x}l_{1,t} U^2{\rm d}t-\big((x+\varepsilon)^\alpha l_{1,x}\big)_x(x+\varepsilon)^\alpha U U_x{\rm d}t,\\
K_3=&-\frac{1}{2}(x+\varepsilon)^\alpha U_x^2+\frac{1}{2}(x+\varepsilon)^\alpha l_{1,x}^2 U^2-\frac{1}{2}l_{1,t}U^2,\\
K_4=&\frac{1}{2}(x+\varepsilon)^\alpha\big({\rm d}U_x\big)^2-\frac{1}{2}(x+\varepsilon)^\alpha l_{1,x}^2({\rm d}U)^2+\frac{1}{2}l_{1,t}({\rm d}U)^2.
\end{align*}

{\em Step 2.\ Carleman estimate for degenerate part.}

In this step, we will prove the Carleman estimate for degenerate part $\big(0, x^{(2)}_1\big)\times (0,T)$:
\begin{align}\label{3.12}
&\mathbb E\int_{0}^T\int_{0}^{x_1^{(2)}} s^3(x+\varepsilon)^{2\alpha+3\beta-4}\xi^3|u^\varepsilon|^2e^{2s\varphi}{\rm d}x {\rm d}t+\mathbb E\int_{0}^T\int_{0}^{x_1^{(2)}} s  (x+\varepsilon)^{2\alpha+\beta-2}\xi |u^\varepsilon_x|^2e^{2s\varphi}{\rm d}x{\rm d}t\nonumber\\
\leq &C\mathbb E\int_{0}^T\int_0^{x^{(2)}_2} |f_1|^2e^{2s\varphi}{\rm d}x{\rm d}t+C(\lambda)\mathbb E \int_{0}^T\int_0^{x^{(2)}_2} s^2\xi^2 |F_1|^2e^{2s\varphi}{\rm d}x{\rm d}t\nonumber\\
&+C\mathbb E \int_{\omega^{(2)}_T}(|u^\varepsilon|^2+|u^\varepsilon_x|^2)e^{2s\varphi}{\rm d}x{\rm d}t+C(\lambda)\mathbb E\int_0^T s^2\left[(x+\varepsilon)^{\alpha+\beta-1}\xi^3 |u^\varepsilon|^2e^{2s\varphi}\right]_{x=0}{\rm d}t.\end{align}

By using
\begin{align}\label{3.13}\left\{\begin{array}{l}
l_{1,x}=\beta s\lambda(x+\varepsilon)^{\beta-1}\phi_1\xi,\quad l_{1,xt}=\beta s\lambda(x+\varepsilon)^{\beta-1}\phi_1\xi_t,\\
 l_{1,xx}=s\left(\beta^2\lambda^2(x+\varepsilon)^{2\beta-2}+\beta(\beta-1)\lambda(x+\varepsilon)^{\beta-2}\right)\phi_1\xi,\\
 |l_{1,t}|=|s\psi_1\xi_t|\leq C(\lambda)s\xi^\frac{3}{2},\quad |l_{1,tt}|=|s\psi_1\xi_{tt}|\leq C(\lambda)s\xi^{2},
\end{array}
\right.
\end{align}
we obtain
\begin{align}\label{3.14}
K_1{\rm d}t\geq &\left(\alpha+2\beta-2\right)\beta^3 s^3\lambda^3(x+\varepsilon)^{2\alpha+3\beta-4}\phi_1^3\xi^3|U|^2 {\rm d}t\nonumber\\
&+(\alpha+2\beta-2)\beta s\lambda(x+\varepsilon)^{2\alpha+\beta-2}\phi_1\xi |U_x|^2{\rm d}t+X_1{\rm d}t+X_2{\rm d}t+X_3{\rm d}t,
\end{align}
where
\begin{align*}
X_1=&-2s^2(x+\varepsilon)^{\alpha}\varphi_{1,x}\varphi_{1,xt}|U|^2,\\
X_2=&\frac{1}{2}s\varphi_{1,tt}|U|^2,\\
X_3=& s(x+\varepsilon)^{\alpha}\left[(x+\varepsilon)^{\alpha}\varphi_{1,x}\right]_{xx}U U_x.
\end{align*}
Now we estimate $X_1$, $X_2$ and $X_3$.
Obviously, by (\ref{3.13}) we have
\begin{align}\label{3.15}
X_1\geq & -Cs^2\lambda^2(x+\varepsilon)^{\alpha+2\beta-2}\phi_1^2\xi^\frac{5}{2}|U|^2\geq-C(\lambda)s^2 (x+\varepsilon)^{2\alpha+3\beta-4}\phi_1^3\xi^3|U|^2,
\end{align}
due to $\beta\leq 2-\alpha$. Obviously,
\begin{align}\label{3.16}
X_2\geq-C(\lambda)s\xi^2|U|^2.
\end{align}
For $X_3$,  we have
\begin{align}\label{3.17}
X_3
 \geq & -\Big[ C_{\alpha,\beta}^{(1)} s\lambda  (x+\varepsilon)^{2\alpha+\beta-3}+Cs\lambda^2(x+\varepsilon)^{2\alpha+2\beta-3} \nonumber\\
 &+Cs\lambda^3(x+\varepsilon)^{2\alpha+3\beta-3} \Big]\phi_1\xi |U||U_x|\nonumber\\
\geq &  -C_{\alpha,\beta}^{(1)} s\lambda  (x+\varepsilon)^{2\alpha+\beta-3}\phi_1\xi |U||U_x|-Cs\lambda^3(x+\varepsilon)^{2\alpha+2\beta-3}\phi_1 \xi |U||U_x|\nonumber\\
\geq &-C_{\alpha,\beta}^{(1)} s\lambda  (x+\varepsilon)^{2\alpha+\beta-3}\phi_1\xi |U||U_x|-Cs(x+\varepsilon)^{2\alpha+\beta-2}\phi_1\xi|U_x|^2\nonumber\\
&-Cs\lambda^6 (x+\epsilon)^{2\alpha+3\beta-4}\phi_1\xi|U|^2.
\end{align}
with \begin{align*}C_{\alpha,\beta}^{(1)}=\beta (\alpha+\beta-1)(2-\alpha-\beta)\geq 0. \end{align*}
Then substituting (\ref{3.14})-(\ref{3.17}) into (\ref{1-3.11}), we find that
\begin{align}\label{3.18}
&\theta_1 I_2\left[{\rm d}u+\big((x+\varepsilon)^\alpha u_x\big)_x {\rm d}t\right]+C(\lambda)s\xi^2|U|^2{\rm d}t+C_{\alpha,\beta}^{(1)} s\lambda  (x+\varepsilon)^{2\alpha+\beta-3}\phi_1\xi |U||U_x|{\rm d}t\nonumber\\
\geq& |I_2|^2{\rm d}t+ \big[\left(\alpha+2\beta-2\right)\beta^3 s^3\lambda^3-C(\lambda)s^2\big] (x+\varepsilon)^{2\alpha+3\beta-4}\phi_1^3\xi^3|U|^2 {\rm d}t\nonumber\\
&+\big[\left(\alpha+2\beta-2\right)\beta s \lambda -Cs \big] (x+\varepsilon)^{2\alpha+\beta-2}\phi_1\xi|U_x|^2{\rm d}t+(K_2)_x\nonumber\\
&+{\rm d}K_3+K_4.
\end{align}
Integrating both side of (\ref{3.18}) on $Q_T$ and taking mathematical expectation, we have
\begin{align}\label{3.19}
&\mathbb E\int_{Q_T}|I_2|^2{\rm d}x{\rm d}t+\mathbb E\int_{Q_T} \big[\left(\alpha+2\beta-2\right)\beta^3 s^3\lambda^3-C(\lambda)s^2\big] (x+\varepsilon)^{2\alpha+3\beta-4}\phi_1^3\xi^3|U|^2{\rm d}x {\rm d}t\nonumber\\
&+\mathbb E\int_{Q_T} \left[\left(\alpha+2\beta-2\right)\beta s \lambda -Cs \right] (x+\varepsilon)^{2\alpha+\beta-2}\phi_1\xi|U_x|^2{\rm d}x{\rm d}t\nonumber\\
\leq& Y_1+Y_2+Y_3-\mathbb E\int_{0}^T\big[K_2\big]_{x=0}^{x=1}-\mathbb E\int_{Q_T}{\rm d}K_3{\rm d}x-\mathbb E\int_{Q_T}K_4{\rm d}x,
\end{align}
where
\begin{align*}
Y_1=&\mathbb E\int_{Q_T}\theta_1 I_2\left[{\rm d}u^\varepsilon+\big((x+\varepsilon)^\alpha u^\varepsilon_x\big)_x{\rm d}t\right]{\rm d}x,\\
Y_2=&C_{\alpha,\beta}^{(1)}\mathbb E \int_{Q_T} s\lambda  (x+\varepsilon)^{2\alpha+\beta-3}\phi_1 \xi |U||U_x|{\rm d}x{\rm d}t,\\
Y_3=&C(\lambda)\mathbb E\int_{Q_T}s\xi^2|U|^2{\rm d}x{\rm d}t.
\end{align*}

 Now we estimate $Y_1,Y_2,Y_3$. For $Y_1$, by noting $\mathbb E\int_{Q_T}\theta_1 I_2 F_1{\rm d}x{\rm d}B(t)=0$, we obtain
\begin{align}\label{3.21}
&Y_1=\mathbb E\int_{Q_T} \theta_1 I_2 (f_1 {\rm d}t+F_1{\rm d}B(t)){\rm d}x\leq \frac{1}{2}\mathbb E\int_{Q_T} |I_2|^2{\rm d}x{\rm d}t+\frac{1}{2}\mathbb E\int_{Q_T}\theta_1^2 |f_1|^2{\rm d}x{\rm d}t.
\end{align}
By using Young's inequality and Lemma 3.2, we have
\begin{align}\label{3.22}
Y_2\leq & \frac{1}{4\epsilon_1} C_{\alpha,\beta}^{(1)} \mathbb E\int_{Q_T} s\lambda(x+\varepsilon)^{-(4-2\alpha-\beta)}\phi_1\xi |U|^2{\rm d}x{\rm d}t\nonumber\\
&+\epsilon_1 C_{\alpha,\beta}^{(1)}\mathbb E\int_{Q_T} s\lambda(x+\varepsilon)^{2\alpha+\beta-2}\phi_1\xi|U_x|^2{\rm d}x{\rm d}t\nonumber\\
\leq &\frac{1}{4\epsilon_1} C_{\alpha,\beta}^{(1)}C_{\alpha,\beta}^{(2)} \mathbb E\int_{Q_T} s\lambda(x+\varepsilon)^{2\alpha+\beta-2}\xi \big|\phi_1^{\frac{1}{2}} U_x+(\phi_1^{\frac{1}{2}})_x U\big|^2{\rm d}x{\rm d}t\nonumber\\
&+\epsilon_1 C_{\alpha,\beta}^{(1)}\mathbb E\int_{Q_T} s\lambda(x+\varepsilon)^{2\alpha+\beta-2}\phi_1\xi|U_x|^2{\rm d}x{\rm d}t\nonumber\\
\leq & \left(\epsilon_1 C_{\alpha,\beta}^{(1)}+\frac{1}{4\epsilon_1} C_{\alpha,\beta}^{(1)}C_{\alpha,\beta}^{(2)}+\epsilon_2\right)\mathbb E\int_{Q_T} s\lambda(x+\varepsilon)^{2\alpha+\beta-2} \phi_1\xi| U_x|^2{\rm d}x{\rm d}t\nonumber\\
&+C(\epsilon_1,\epsilon_2)\mathbb E\int_{Q_T} s\lambda^3(x+\varepsilon)^{2\alpha+3\beta-4} \phi_1\xi |U|^2{\rm d}x{\rm d}t,
\end{align}
with
\begin{align*}
C_{\alpha,\beta}^{(2)}=\frac{4}{(3-2\alpha-\beta)^2}.
\end{align*}
Similarly,
\begin{align}\label{3.23}
Y_3\leq &C \mathbb E\int_{Q_T} (x+\varepsilon)^{-(4-2\alpha-\beta)}\xi |U|^2{\rm d}x{\rm d}t+C(\lambda)\mathbb E\int_{Q_T} s^2(x+\varepsilon)^{4-2\alpha-\beta} \xi^3 |U|^2{\rm d}x{\rm d}t\nonumber\\
\leq & C\mathbb E\int_{Q_T} (x+\varepsilon)^{2\alpha+\beta-2}\xi |U_x|^2{\rm d}x{\rm d}t+C(\lambda)\mathbb E\int_{Q_T} s^2(x+\varepsilon)^{2\alpha+3\beta-4} \xi^3 |U|^2{\rm d}x{\rm d}t.
\end{align}

From (\ref{3.19})-(\ref{3.23}), it follows that
\begin{align}\label{3.26}
&\mathbb E\int_{Q_T} \big[(\alpha+2\beta-2)\beta^3 s^3\lambda^3-C(\lambda)s^2-C(\epsilon_1,\epsilon_2)s\lambda^3\big](x+\varepsilon)^{2\alpha+3\beta-4}\phi_1^3\xi^3|U|^2{\rm d}x {\rm d}t\nonumber\\
&+\mathbb E\int_{Q_T} \big(C_{\alpha,\beta}^{(3)} s \lambda -Cs -C\big) (x+\varepsilon)^{2\alpha+\beta-2}\phi_1\xi|U_x|^2{\rm d}x{\rm d}t\nonumber\\
\leq &C\mathbb E\int_{Q_T}\theta_1^2 |f_1|^2{\rm d}x{\rm d}t-\mathbb E\int_0^T\big[K_2\big]_{x=0}^{x=1}-\mathbb E\int_{Q_T}{\rm d}K_3{\rm d}x-\mathbb E\int_{Q_T}K_4{\rm d}x,
\end{align}
where
\begin{align*}
C_{\alpha,\beta}^{(3)} =(\alpha+2\beta-2)\beta-\epsilon_1 C_{\alpha,\beta}^{(1)}-\frac{1}{4\epsilon_1} C_{\alpha,\beta}^{(1)}C_{\alpha,\beta}^{(2)}-\epsilon_2.
\end{align*}
By using (\ref{3.2}),  we can prove for all $\alpha\in (0,2)$ that
\begin{align}\label{3.27}
&\epsilon_1 C_{\alpha,\beta}^{(1)}+\frac{1}{4\epsilon_1} C_{\alpha,\beta}^{(1)}C_{\alpha,\beta}^{(2)}<(\alpha+2\beta-2)\beta.
\end{align}
 We first fix $\epsilon_1=\frac{1}{|3-2\alpha-\beta|}$. For $\alpha\in (0,1)$, (\ref{3.27}) can be simplifies as
 \begin{align*}
 \alpha-\alpha\beta+2\beta-2>0,
 \end{align*}
 which holds for $\beta>1$. For $\alpha\in (1,2)$, since $\beta>3-2\alpha$, (\ref{3.27}) is equivalent to
 \begin{align*}
 4\beta^2+(9\alpha-14)\beta+4\alpha^2-13\alpha +10>0,
 \end{align*}
 which holds for $\beta>\frac{14-9\alpha+\sqrt{17\alpha^2-44\alpha+36}}{8}$.
 Moreover, when $\alpha=1$, we easily see $C^{(1)}_{\alpha,\beta}=0$ and then  (\ref{3.27}). Therefore, (\ref{3.27}) holds for all $\alpha\in (0,2)$, if $\beta$ satisfies (\ref{3.2}).  Further for sufficiently small $\epsilon_2$ we have $C_{\alpha,\beta}^{(3)}>0$.   Consequently, there exist $\lambda_1$ and $s_1$ such that for all $\lambda>\lambda_1$ and $s>s_1$, it holds that
 \begin{align}\label{3.28}
&\mathbb E\int_{Q_T} s^3\lambda^3(x+\varepsilon)^{2\alpha+3\beta-4}\phi_1^3\xi^3\theta_1^2|u^\varepsilon|^2{\rm d}x {\rm d}t+\mathbb E\int_{Q_T} s \lambda(x+\varepsilon)^{2\alpha+\beta-2}\phi_1\xi\theta_1^2|u^\varepsilon_x|^2{\rm d}x{\rm d}t\nonumber\\
\leq &C\mathbb E\int_{Q_T}\theta_1^2 |f_1|^2{\rm d}x{\rm d}t-C\mathbb E\int_0^T\big[K_2\big]_{x=0}^{x=1}-C\mathbb E\int_{Q_T}{\rm d}K_3{\rm d}x-C\mathbb E\int_{Q_T}K_4{\rm d}x.
\end{align}

Now we deal with the boundary term of $K_2$.  For $\alpha\in (0,1)$, by using (\ref{1-3.8}) we have
\begin{align}\label{3.29}
-\mathbb E\int_0^T\big[K_2\big]_{x=0}^{x=1}
= &\mathbb E\int_0^T \big[ (x+\varepsilon)^{2\alpha} l_{1,x} |U_x|^2\big]_{x=0}^{x=1}{\rm d}t\leq  C\mathbb E\int_0^T s\lambda\left[\phi_1\xi \theta_1^2|u^\varepsilon_x|^2\right]_{x=1} {\rm d}t.
\end{align}
Similarly, for $\alpha\in [1,2)$ we have
\begin{align}\label{3.30}
-\mathbb E\int_0^T\big[K_2\big]_{x=0}^{x=1}= &\mathbb E\int_0^T\left[(x+\varepsilon)^\alpha U_x{\rm d}U\right]_{x=0}+\mathbb E\int_0^T\left[(x+\varepsilon)^{2\alpha}l_{1,x} |U_x|^2\right]_{x=0}^{x=1}{\rm d}t\nonumber\\
&-\mathbb E\int_0^T\left[(x+\varepsilon)^{2\alpha}l_{1,x}^3 |U|^2\right]_{x=0}{\rm d}t+\mathbb E\int_0^T\left[(x+\varepsilon)^\alpha l_{1,x}l_{1,t}|U|^2\right]_{x=0}{\rm d}t\nonumber\\
&-\mathbb E\int_0^T\left[\big((x+\varepsilon)^\alpha l_{1,x}\big)_{x}(x+\varepsilon)^\alpha U U_x\right]_{x=0}{\rm d}t\nonumber\\
\leq &\mathbb E\int_0^T\left[(x+\varepsilon)^\alpha U_x{\rm d}U\right]_{x=0}+C\mathbb E\int_0^T s\lambda\left[\phi_1\xi \theta_1^2|u^\varepsilon_x|^2\right]_{x=1} {\rm d}t\nonumber\\
&+C(\lambda)\mathbb E\int_0^T s^2\left[(x+\varepsilon)^{\alpha+\beta-1}\phi_1\xi^3 \theta_1^2 |u^\varepsilon|^2\right]_{x=0}{\rm d}t.
\end{align}
By using It\^{o} formula and (\ref{1-3.8}) again,  we have
\begin{align}\label{3.31}
&\mathbb E\int_0^T\left[(x+\varepsilon)^\alpha U_x{\rm d}U\right]_{x=0}=\mathbb E\int_0^T\left[(x+\varepsilon)^\alpha l_{1,x} U {\rm d}U\right]_{x=0}\nonumber\\
=&\frac{1}{2}\mathbb E\int_0^T\beta s\lambda{\rm d}\left[(x+\varepsilon)^{\alpha+\beta-1}\phi_1\xi |U|^2\right]_{x=0}-\frac{1}{2}\mathbb E\int_0^T\beta s\lambda\left[(x+\varepsilon)^{\alpha+\beta-1}\phi_1\xi_t |U|^2\right]_{x=0}{\rm d}t\nonumber\\
&-\frac{1}{2}\mathbb E\int_0^T \beta s\lambda\left[(x+\varepsilon)^{\alpha+\beta-1}\phi_1\xi ({\rm d}U)^2\right]_{x=0}\nonumber\\
\leq & C(\lambda)\mathbb E\int_0^T s\left[(x+\varepsilon)^{\alpha+\beta-1}\phi_1\xi^2\theta_1^2 |u^\varepsilon|^2\right]_{x=0}{\rm d}t.
\end{align}
Therefore, combining (\ref{3.29})-(\ref{3.31}), we obtain for all $\alpha\in (0,2)$ that
\begin{align}\label{3.32}
&-\mathbb E\int_0^T\big[K_2\big]_{x=0}^{x=1}\nonumber\\
\leq &C\mathbb E\int_0^T s\lambda\left[\phi_1\xi \theta_1^2|u^\varepsilon_x|^2\right]_{x=1} {\rm d}t+C(\lambda)\mathbb E\int_0^T s^2\left[(x+\varepsilon)^{\alpha+\beta-1}\phi_1\xi^3 \theta_1^2 |u^\varepsilon|^2\right]_{x=0}{\rm d}t.
\end{align}
By using $U(x,0)=U(x,T)=0$, we have
\begin{align}\label{3.24}
-\mathbb E\int_{Q_T}{\rm d}K_3{\rm d}x=0.
\end{align}
Moreover, by (\ref{3.13}), $({\rm d}U)^2=\theta_1^2|F_1|^2{\rm d}t$ and $\beta>1-\frac{\alpha}{2}$ for all $\alpha\in (0,2)$,  we have the following estimate:
\begin{align}\label{3.25}
-\mathbb E\int_{Q_T}K_4{\rm d}x\leq& C\mathbb E\int_{Q_T} s^2\lambda^2(x+\varepsilon)^{\alpha+2\beta-2}\phi_1^2\xi^2({\rm d}U)^2{\rm d}x+C(\lambda)\mathbb E\int_{Q_T}s\xi^\frac{3}{2} ({\rm d}U)^2{\rm d}x\nonumber\\
\leq& C(\lambda)\mathbb E\int_{Q_T} s^2\phi_1^2\xi^2\theta_1^2|F_1|^2{\rm d}x{\rm d}t.
\end{align}
 Then substituting (\ref{3.24}) and (\ref{3.25}) into (\ref{3.28}) yields
 \begin{align}\label{1-3.35}
&\mathbb E\int_{Q_T} s^3\lambda^3(x+\varepsilon)^{2\alpha+3\beta-4}\phi_1^3\xi^3\theta_1^2|u^\varepsilon|^2{\rm d}x {\rm d}t+\mathbb E\int_{Q_T} s \lambda(x+\varepsilon)^{2\alpha+\beta-2}\phi_1\xi\theta_1^2|u^\varepsilon_x|^2{\rm d}x{\rm d}t\nonumber\\
\leq &C\mathbb E\int_{Q_T}\theta_1^2 |f_1|^2{\rm d}x{\rm d}t+C(\lambda)\mathbb E \int_{Q_T} s^2\phi_1^2\xi^2\theta_1^2 |F_1|^2{\rm d}x{\rm d}t\nonumber\\
&+C(\lambda)\mathbb E\int_0^T s^2\left[(x+\varepsilon)^{\alpha+\beta-1}\phi_1\xi^3\theta_1^2 |u|^2\right]_{x=0}{\rm d}t+C\mathbb E\int_0^T s\lambda\left[\phi_1\xi \theta_1^2|u^\varepsilon_x|^2\right]_{x=1} {\rm d}t.
\end{align}

Next, we eliminate the boundary term on $x=1$. We consider the following stochastic parabolic equation of $\tilde u^\varepsilon=\chi u^\varepsilon$:
\begin{align}\label{3.33}\left\{
\begin{array}{ll}
{\rm d}\tilde u^\varepsilon+\left((x+\varepsilon)^\alpha \tilde u^\varepsilon_x\right)_x{\rm d}t=\tilde f_1{\rm d}t+\tilde F_1{\rm d}B(t),&(x,t)\in Q_T,\\
 \tilde u^\varepsilon(1,t)=0,& t\in (0,T),\\
 {\rm and}\ \left\{\begin{array}{ll}
 \tilde u^\varepsilon(0,t)=0 & {\rm for}\ \alpha\in (0,1),\\
 ((x+\varepsilon)^\alpha \tilde u^\varepsilon_x)(0,t)=0 & {\rm for}\ \alpha\in [1,2),
 \end{array}\right. &t\in (0,T).
 \end{array}\right.
\end{align}
where
\begin{align*}
\tilde f_1=\left((x+\varepsilon)^{\alpha}\chi_x u^\varepsilon\right)_x+(x+\varepsilon)^{\alpha}\chi_x u^\varepsilon_x+\chi f_1,\quad \tilde F_1=\chi F_1.
\end{align*}
Applying (\ref{1-3.35}) to $\tilde u$ and using the definition of $\chi$, we find that
\begin{align}\label{3.35}
&\mathbb E\int_{0}^T\int_{0}^{x_1^{(2)}} s^3(x+\varepsilon)^{2\alpha+3\beta-4}\xi^3\theta_1^2|u^\varepsilon|^2{\rm d}x {\rm d}t+\mathbb E\int_{0}^T\int_{0}^{x_1^{(2)}} s (x+\varepsilon)^{2\alpha+\beta-2}\xi\theta_1^2|u^\varepsilon_x|^2{\rm d}x{\rm d}t\nonumber\\
\leq &C\mathbb E\int_{Q_T}\chi^2\theta_1^2 |f_1|^2{\rm d}x{\rm d}t+C(\lambda)\mathbb E \int_{Q_T} s^2\chi^2\xi^2\theta_1^2 |F_1|^2{\rm d}x{\rm d}t\nonumber\\
&+C\mathbb E\int_{\omega^{(2)}_T}\theta_1^2 (|u^\varepsilon|^2+|u^\varepsilon_x|^2){\rm d}x{\rm d}t+C(\lambda)\mathbb E\int_0^T s^2\left[(x+\varepsilon)^{\alpha+\beta-1}\xi^3\theta_1^2 |u^\varepsilon|^2\right]_{x=0}{\rm d}t.
\end{align}
Together with
$\varphi_1=\varphi$ for $x\in (0,x_2^{(2)})$, we deduce (\ref{3.12}) from (\ref{3.35}).

{\em Step 3. Carleman estimate for nondegenerate part.}

Now we derive the Carleman estimate for nondegenerate part
$(x_2^{(2)},1)\times (0,T)$:
\begin{align}\label{3.36}
&\mathbb E\int_{0}^T\int_{x_2^{(2)}}^1 s^3\xi^3(x+\varepsilon)^{2\alpha+3\beta-4}|u^\varepsilon|^2e^{2s\varphi}{\rm d}x{\rm d}t+\mathbb E \int_{0}^T\int_{x_2^{(2)}}^1 s\xi (x+\varepsilon)^{2\alpha+\beta-2} |u^\varepsilon_x|^2e^{2s\varphi}{\rm d}x{\rm d}t\nonumber\\
\leq &C\mathbb E\int_{0}^T\int_{x^{(2)}_1}^1| f_1|^2e^{2s\varphi}{\rm
d}x{\rm
d}t+C(\lambda)\mathbb E\int_{0}^T\int_{x^{(2)}_1}^1s^2\xi^2|F_1|^2e^{2s\varphi}{\rm d}x{\rm
d}t\nonumber\\
&+C(\lambda)\mathbb E\int_{\omega^{(1)}_T} \left(|u^\varepsilon_x|^2+s^3\xi^3|u^\varepsilon|^2\right)e^{2s\varphi}{\rm
d}x{\rm d}t.
\end{align}
 To do this, letting $\overline u^\varepsilon=(1-\chi)u^\varepsilon$, then
we have
\begin{align}\label{3.37}\left\{
\begin{array}{ll}
{\rm d}\overline u^\varepsilon+\left((x+\varepsilon)^\alpha \overline u^\varepsilon_x\right)_x{\rm d}t=\overline{g}_1{\rm d}t+\overline{G}_1{\rm d}B(t),&(x,t)\in Q_T,\\
 \overline u^\varepsilon(1,t)=0,& t\in (0,T),\\
 {\rm and}\ \left\{\begin{array}{ll}
 \overline u^\varepsilon(0,t)=0 & {\rm for}\ \alpha\in (0,1),\\
 ((x+\varepsilon)^\alpha \overline u^\varepsilon_x)(0,t)=0 & {\rm for}\ \alpha\in [1,2),
 \end{array}\right. &t\in (0,T).
 \end{array}\right.
\end{align}
where
\begin{align*}
\overline{g}_1=(1-\chi)f_1-((x+\varepsilon)^{\alpha}\chi_x
u^\varepsilon)_x-(x+\varepsilon)^\alpha \chi_x u^\varepsilon_x,\quad \overline{F_1}=(1-\chi)F_1.
\end{align*} By the classic Carleman estimate for  stochastic nondegenerate parabolic equation, e.g. [\ref{Liu2014ECOCV}] or [\ref{Yan2015IP}],  we have
\begin{align}\label{3.38}
&\mathbb E\int_{Q_T}s^3\lambda^3\xi^3|\overline u^\varepsilon|^2e^{2s\varphi_2}{\rm d}x{\rm d}t+\mathbb E \int_{Q_T} s\lambda\xi|\overline u^\varepsilon_x|^2e^{2s\varphi_2}{\rm d}x{\rm d}t\nonumber\\
\leq& C\mathbb E\int_{Q_T}|\overline f_1|^2e^{2s\varphi_2}{\rm
d}x{\rm
d}t+C\mathbb E\int_{Q_T}
s^2\lambda^2\xi^2|\overline F_1|^2e^{2s\varphi_2}{\rm d}x{\rm
d}t+C\mathbb E\int_{\omega^{(1)}_T} s^3\lambda^3\xi^3|\overline u^\varepsilon|^2e^{2s\varphi_2}{\rm
d}x{\rm d}t\nonumber\\
\leq &C\mathbb E\int_{Q_T}(1-\chi)^2| f_1|^2e^{2s\varphi_2}{\rm
d}x{\rm
d}t+C\mathbb E\int_0^T\int_{x^{(2)}_1}^{x_2^{(2)}}(|u^\varepsilon_x|^2+|u^\varepsilon|^2)e^{2s\varphi_2}{\rm d}x{\rm d}t\nonumber\\
&+C\mathbb E\int_{Q_T}
s^2\lambda^2(1-\chi)^2\xi^2|F_1|^2e^{2s\varphi_2}{\rm d}x{\rm
d}t+C\mathbb E\int_0^T\int_{x^{(2)}_1}^{x^{(1)}_2} s^3\lambda^3\xi^3|u^\varepsilon|^2e^{2s\varphi_2}{\rm
d}x{\rm d}t.
\end{align}
Since $\varphi_2=\varphi$ for
$x\in (x_1^{(2)},1)$ and $\min\{(x+\varepsilon)^{2\alpha+3\beta-4},(x+\varepsilon)^{2\alpha+\beta-2}\}\geq C>0$ for
$x\in (x_1^{(2)},1)$,  together with (\ref{3.38}),  we obtain (\ref{3.36}).

{\em Step 4. End of the proof.}

Combining  (\ref{3.12}) and (\ref{3.36}) and adding to both sides of
the inequality the term
\begin{align*}
\mathbb E\int_0^T\int_{x_1^{(2)}}^{x_2^{(2)}}s^3\xi^3(x+\varepsilon)^{2\alpha+3\beta-4}|u^\varepsilon|^2e^{2s\varphi}{\rm d}x{\rm d}t+\mathbb E\int_0^T\int_{x_1^{(2)}}^{x_2^{(2)}}s\xi(x+\varepsilon)^{2\alpha+\beta-2}
|u^\varepsilon_x|^2e^{2s\varphi}{\rm d}x{\rm d}t,
\end{align*}
 we obtain
\begin{align}
&\mathbb E\int_{Q_T} s^3\lambda^3\xi^3(x+\varepsilon)^{2\alpha+3\beta-4}|u^\varepsilon|^2e^{2s\varphi}{\rm d}x {\rm d}t+\mathbb E\int_{Q_T} s \lambda\xi (x+\varepsilon)^{2\alpha+\beta-2}|u^\varepsilon_x|^2 e^{2s\varphi}{\rm d}x{\rm d}t\nonumber\\
\leq & C\mathbb E\int_{Q_T}|f_1|^2e^{2s\varphi}{\rm d}x{\rm
d}t+C(\lambda)\mathbb E\int_{Q_T}
s^2\xi^2|F_1|^2e^{2s\varphi}{\rm d}x{\rm
d}t+C(\lambda)\mathbb E\int_{\omega^{(1)}_T} s^3\xi^3|u^\varepsilon|^2e^{2s\varphi}{\rm
d}x{\rm d}t\nonumber\\
&+C(\lambda)\mathbb E\int_{\omega^{(1)}_T} s\xi|u^\varepsilon_x|^2e^{2s\varphi}{\rm
d}x{\rm d}t+C(\lambda)\mathbb E\int_0^T s^2\left[(x+\varepsilon)^{\alpha+\beta-1}\xi^3 |u^\varepsilon|^2e^{2s\varphi}\right]_{x=0}{\rm d}t.
\end{align}

Finally, by the Cacciopoli inequality (\ref{3.7}), we obtain (\ref{3.3}). This completes the proof of Theorem 3.1.\hfill $\Box$

\subsection{Carleman estimate for forward stochastic degenerate equation}

In this subsection, we will introduce a regular weight function into a new Carleman estimate for the backward stochastic degenerate equation, in which the random source and the initial data are put on the left-hand side. This allows us to prove the stability for our inverse  random source problem.

We set
\begin{align*}
\varrho_i(x,t)=\eta_i(x)-(\lambda-t)^2+\lambda^2,\quad \Phi_i(x,t)=e^{\lambda\varrho_i(x,t)},\quad i=1,2,
\end{align*}
where $\eta_i$ $(i=1,2)$ are same as the ones in section 3.1.  We introduce regular weight function
\begin{align*}
&\Phi(x,t)=\chi(x)\Phi_1(x,t)+(1-\chi(x))\Phi_2(x,t),\quad (x,t)\in Q_T.
\end{align*}
So that, similar to (\ref{3.1}) we also have
\begin{align}\label{3.42}
\Phi(x,t)=\left\{\begin{array}{ll}\Phi_1(x,t),&(x,t)\in (0,x_1^{(2)})\times (0,T),\\
\Phi_1(x,t)=\Phi_2(x,t),&(x,t)\in (x_1^{(2)},x_2^{(2)})\times (0,T),\\
\Phi_2(x,t),&(x,t)\in (x_2^{(2)},1)\times (0,T).\\
\end{array}\right.
\end{align}

\begin{thm} Let $\alpha\in \left(0,2\right)$,  $f_2\in L^2_\mathcal F(0,T;L^2(I))$, $F_2\in L^2_\mathcal F(0,T; H^1(I))$ and $\beta$ such that (\ref{3.2}). Then for any $\varepsilon\in (0,1)$, there exist positive constants $\lambda_2=\lambda_2(\omega,I,T,\alpha)$, $s_2=s_2(\omega,I,T,\alpha,\lambda)$ and $C_5=C_5(\omega,I,T$, $\alpha)$, $C_6=C_6(\omega,I,T,\alpha,\lambda)$
such that
\begin{align}\label{3.43}
&\mathbb E\int_{Q_T}s\lambda\Phi|F_2|^2e^{2s\Phi}{\rm d}x{\rm d}t+\mathbb E\int_{Q_T} s^3\lambda^3\Phi^3(x+\varepsilon)^{2\alpha+3\beta-4} |v^\varepsilon|^2e^{2s\Phi}{\rm d}x {\rm d}t\nonumber\\
&+\mathbb E\int_{Q_T} s\lambda\Phi(x+\varepsilon)^{2\alpha+\beta-2}|v^\varepsilon_x|^2e^{2s\Phi}{\rm d}x{\rm d}t\nonumber\\
\leq &C_5\mathbb E\int_{Q_T} |f_2|^2e^{2s\Phi}{\rm d}x{\rm d}t+C_5\mathbb E\int_{Q_T} s\Phi|F_{2,x}|^2e^{2s\Phi}{\rm d}x{\rm d}t\nonumber\\
&+C_6(\lambda)\mathbb E\int_{\omega_T}s^3\Phi^3|v^\varepsilon|^2e^{2s\Phi}{\rm d}x{\rm d}t+C_6(\lambda)s^2 e^{C(\lambda)s}\|v^\varepsilon(\cdot, T)\|_{L^2(\Omega,\mathcal F_T,P;L^2(I))}^2\nonumber\\
&+C_6(\lambda)\mathbb E\int_0^T s^2\left[(x+\varepsilon)^{\alpha+\beta-1} (|v^\varepsilon|^2+|F_2|^2)e^{2s\Phi}\right]_{x=0}{\rm d}t
\end{align}
for all $\lambda\geq \lambda_2$, $s\geq s_2$ and all $v^\varepsilon\in \mathcal H^1$ satisfying
\begin{align}\left\{
\begin{array}{ll}
{\rm d}v^\varepsilon-\left((x+\varepsilon)^\alpha v^\varepsilon_x\right)_x{\rm d}t=f_2{\rm d}t+F_2{\rm d}B(t),&(x,t)\in Q_T,\\
 v^\varepsilon(1,t)=0,& t\in (0,T),\\
 {\rm and}\ \left\{\begin{array}{ll}
 v^\varepsilon(0,t)=0 & {\rm for}\ \alpha\in (0,1),\\
 ((x+\varepsilon)^\alpha v^\varepsilon_x)(0,t)=0 & {\rm for}\ \alpha\in [1,2),
 \end{array}\right. &t\in (0,T),\\
 v^\varepsilon(x,0)=0,&x\in I.
 \end{array}\right.
\end{align}
\end{thm}

{\noindent\bf Remark 3.3.}\ The second large parameter $\lambda$ in the proof of null controllability could be omitted. However,  in inverse random source problem it  plays a very important role.

\vspace{2mm}

Based on Theorem 3.5, letting $\beta=2-\alpha$ and $\varepsilon\rightarrow 0$, we could drop the boundary term in (\ref{3.43}) as in Theorem 3.2. Then we obtain the following result:

\begin{thm} Let $\alpha\in \left(0,2\right)$,  $f_2\in L^2_\mathcal F(0,T;L^2(I))$, $F_2\in L^2_\mathcal F(0,T; H^1(I))$ and $\beta$ such that (\ref{3.2}). Then for any $\varepsilon\in (0,1)$, there exist positive constants $\lambda_2=\lambda_2(\omega,I,T,\alpha)$, $s_2=s_2(\omega,I,T,\alpha,\lambda)$ and $C_5=C_5(\omega,I,T$, $\alpha)$, $C_6=C_6(\omega,I,T,\alpha,\lambda)$
such that
\begin{align}\label{3.45}
&\mathbb E\int_{Q_T}s\lambda\Phi|F_2|^2e^{2s\Phi}{\rm d}x{\rm d}t+\mathbb E\int_{Q_T} s^3\lambda^3\Phi^3 x^{2-\alpha} |v|^2e^{2s\Phi}{\rm d}x {\rm d}t\nonumber\\
&+\mathbb E\int_{Q_T} s\lambda\Phi x^{\alpha}|v_x|^2e^{2s\Phi}{\rm d}x{\rm d}t\nonumber\\
\leq &C_5\mathbb E\int_{Q_T} |f_2|^2e^{2s\Phi}{\rm d}x{\rm d}t+C_5\mathbb E\int_{Q_T} s\Phi|F_{2,x}|^2e^{2s\Phi}{\rm d}x{\rm d}t\nonumber\\
&+C_6(\lambda)\mathbb E\int_{\omega_T}s^3\Phi^3|v|^2e^{2s\Phi}{\rm d}x{\rm d}t+C_6(\lambda)s^2 e^{C(\lambda)s}\|v(\cdot, T)\|_{L^2(\Omega,\mathcal F_T,P;L^2(I))}^2
\end{align}
for all $\lambda\geq \lambda_2$, $s\geq s_2$ and all $v\in \mathcal H_\alpha^1$ satisfying
\begin{align}\left\{
\begin{array}{ll}
{\rm d}v-\left(x^\alpha v_x\right)_x{\rm d}t=f_2{\rm d}t+F_2{\rm d}B(t),&(x,t)\in Q_T,\\
 v(1,t)=0,& t\in (0,T),\\
 {\rm and}\ \left\{\begin{array}{ll}
 v(0,t)=0 & {\rm for}\ \alpha\in (0,1),\\
 \big(x^\alpha v_x\big)(0,t)=0 & {\rm for}\ \alpha\in [1,2),
 \end{array}\right. &t\in (0,T),\\
 v(x,0)=0,&x\in I.
 \end{array}\right.
\end{align}
\end{thm}

Now we prove Theorem 3.5.

\vspace{2mm}

\noindent{\bf Proof of Theorem 3.5.}\   Let $L_1=s\Phi_1$, $\Theta_1=e^{L_1}$ and $V=\Theta_1 v^\varepsilon$. Then we have
\begin{align}
\Theta_1 \left[{\rm d}v^\varepsilon-\big((x+\varepsilon)^\alpha v^\varepsilon_x\big)_x{\rm d}t\right]=J_1+J_2{\rm d}t
\end{align}
with
\begin{align}\label{}\left\{\begin{array}{ll}
V(1,t)=0,& t\in (0,T),\\
{\rm and}\ \left\{\begin{array}{ll}
V(0,t)=0& {\rm for}\ \alpha\in (0,1),\\
((x+\varepsilon)^\alpha V_x)(0,t)=\big((x+\varepsilon)^{\alpha}l_{1,x}U\big)(0,t)& {\rm for}\ \alpha\in [1,2),
\end{array}
\right. &t\in (0,T)\\
V(x,0)=0,& x\in \overline I,\end{array}\right.
\end{align}
where
\begin{align*}
&J_1={\rm d}V+2(x+\varepsilon)^\alpha L_{1,x}V_x{\rm d}t+\big((x+\varepsilon)^\alpha L_{1,x}\big)_x V{\rm d}t,\\
&J_2= -\big((x+\varepsilon)^\alpha V_{x}\big)_x-(x+\varepsilon)^{\alpha}L_{1,x}^2 V-L_{1,t} V.
\end{align*}
Hence,
\begin{align}
\Theta_1 J_2\left[{\rm d}v^\varepsilon-\big((x+\varepsilon)^\alpha v^\varepsilon_x\big)_x{\rm d}t\right]=J_1J_2+|J_2|^2{\rm d}t.
\end{align}
Then by a similar argument to (\ref{1-3.11}), we have
\begin{align}
&\Theta_1 J_2\big[{\rm d}v^\varepsilon-\big((x+\varepsilon)^\alpha v^\varepsilon_x\big)_x{\rm d}t\big]=|J_2|^2{\rm d}t+R_1{\rm d}t+(R_2)_x+{\rm d}R_3+R_4,
\end{align}
where
\begin{align*}
R_1=&\left[\alpha(x+\varepsilon)^{2\alpha-1}L_{1,x}^3+2(x+\varepsilon)^{2\alpha}L_{1,x}^2 L_{1,xx}\right]V^2\\
&+\big[2(x+\varepsilon)^{2\alpha}L_{1,xx}+\alpha(x+\varepsilon)^{2\alpha-1}L_{1,x}\big]V_x^2+\big[(x+\varepsilon)^\alpha L_{1,x}\big]_{xx} (x+\varepsilon)^\alpha VV_x\nonumber\\
&+(x+\varepsilon)^\alpha L_{1,x}L_{1,xt}V^2+\frac{1}{2}L_{1,tt} V^2+\big[(x+\varepsilon)^\alpha L_{1,x}L_{1,t}\big]_x V^2{\rm d}t\nonumber\\
&-\big[(x+\varepsilon)^\alpha L_{1,x}\big]_x L_{1,t} V^2,\nonumber\\
R_2=&-(x+\varepsilon)^\alpha V_x{\rm d}V-(x+\varepsilon)^{2\alpha} L_{1,x} V_x^2{\rm d}t-(x+\varepsilon)^{2\alpha}L_{1,x}^3 V^2{\rm d}t\\
&-(x+\varepsilon)^\alpha L_{1,x}L_{1,t} V^2{\rm d}t-\big[(x+\varepsilon)^\alpha L_{1,x}\big]_x(x+\varepsilon)^\alpha V V_x{\rm d}t,\\
R_3=&\frac{1}{2}(x+\varepsilon)^\alpha V_x^2-\frac{1}{2}(x+\varepsilon)^\alpha L_{1,x}^2 V^2-\frac{1}{2}L_{1,t}V^2,\\
R_4=&-\frac{1}{2}(x+\varepsilon)^\alpha\big({\rm d}V_x\big)^2+\frac{1}{2}(x+\varepsilon)^\alpha L_{1,x}^2({\rm d}V)^2+\frac{1}{2}L_{1,t}({\rm d}V)^2.
\end{align*}
Noticing that the weight function is regular, we have
\begin{align}\label{}\left\{\begin{array}{l}
L_{1,x}=\beta s\lambda(x+\varepsilon)^{\beta-1}\Phi_1,\quad L_{1,t}=2 s\lambda (\lambda-t)\Phi_1,\\
 L_{1,xt}=2\beta s\lambda^2(x+\varepsilon)^{\beta-1}(\lambda-t)\Phi_1,\\ L_{1,xx}=s\left(\beta^2\lambda^2(x+\varepsilon)^{2\beta-2}+\beta(\beta-1)\lambda(x+\varepsilon)^{\beta-2}\right)\Phi_1,\\
 L_{1,tt}=s\left(4\lambda^2(\lambda-t)^2-2\lambda\right)\Phi_1.
\end{array}
\right.
\end{align}
Then by a similar process to obtain (\ref{3.28}), we can prove that  there exist $\lambda_2$ and $s_2$ such that for all $\lambda>\lambda_2$ and $s>s_2$, it holds that
 \begin{align}\label{3.50}
&\mathbb E\int_{Q_T} s^3\lambda^3(x+\varepsilon)^{2\alpha+3\beta-4}\Phi_1^3\Theta_1^2|v^\varepsilon|^2{\rm d}x {\rm d}t+\mathbb E\int_{Q_T} s \lambda(x+\varepsilon)^{2\alpha+\beta-2}\Phi_1\Theta_1^2|v^\varepsilon_x|^2{\rm d}x{\rm d}t\nonumber\\
\leq &C\mathbb E\int_{Q_T}\Theta_1^2 |f_2|^2{\rm d}x{\rm d}t-\mathbb E\int_0^T\big[R_2\big]_{x=0}^{x=1}-\mathbb E\int_{Q_T}{\rm d}R_3{\rm d}x-\mathbb E\int_{Q_T}R_4{\rm d}x.
\end{align}

Now we analyze the terms of $R_2$, $R_3$ and $R_4$. For the boundary term of $R_2$, noticing that
\begin{align*}
&-\mathbb E\int_0^T\left[(x+\varepsilon)^\alpha V_x{\rm d}V\right]_{x=0}=-\mathbb E\int_0^T\left[(x+\varepsilon)^\alpha L_{1,x} V {\rm d}V\right]_{x=0}\nonumber\\
=&-\frac{1}{2}\mathbb E\int_0^T\beta s\lambda{\rm d}\left[(x+\varepsilon)^{\alpha+\beta-1}\Phi_1|V|^2\right]_{x=0}+\mathbb E\int_0^T\beta s\lambda^2(\lambda-t)\left[(x+\varepsilon)^{\alpha+\beta-1}\Phi_1 |V|^2\right]_{x=0}{\rm d}t\nonumber\\
&+\frac{1}{2}\mathbb E\int_0^T \beta s\lambda\left[(x+\varepsilon)^{\alpha+\beta-1}\Phi_1 ({\rm d}V)^2\right]_{x=0}\nonumber\\
\leq & C(\lambda)\mathbb E\int_0^T s\left[(x+\varepsilon)^{\alpha+\beta-1}\Phi_1\Theta_1^2 |v^\varepsilon|^2\right]_{x=0}{\rm d}t+C\mathbb E\int_0^Ts\lambda\left[(x+\varepsilon)^{\alpha+\beta-1}\Phi_1\Theta_1^2|F_2|^2\right]_{x=0}{\rm d}t,
\end{align*}
similar to (\ref{3.32}), we obtain
\begin{align}\label{3.51}
&-\mathbb E\int_0^T\big[R_2\big]_{x=0}^{x=1}\nonumber\\
\leq & C\mathbb E\int_0^T s\lambda\left[\Phi_1\Theta_1^2|v^\varepsilon_x|^2\right]_{x=1} {\rm d}t+C(\lambda)\mathbb E\int_0^T s^2\left[(x+\varepsilon)^{\alpha+\beta-1}\Phi_1^2 \Theta_1^2 |v^\varepsilon|^2\right]_{x=0}{\rm d}t\nonumber\\
&+C\mathbb E\int_0^Ts\lambda\left[(x+\varepsilon)^{\alpha+\beta-1}\Phi_1\Theta_1^2|F_2|^2\right]_{x=0}{\rm d}t.
\end{align}
By $V(x,0)=0$, $\mathbb P$-a.s.  in $I$,   we have
\begin{align}\label{3.52}
-\mathbb E\int_{Q_T}{\rm d}R_3{\rm d}x=&\mathbb E\int_{I}\left[-\frac{1}{2}(x+\varepsilon)^\alpha V_x^2+\frac{1}{2}(x+\varepsilon)^\alpha L_{1,x}^2 V^2+\frac{1}{2}L_{1,t}V^2\right]_{t=T}{\rm d}x\nonumber\\
\leq & C(\lambda)s^2 e^{C(\lambda)s}\|v^\varepsilon(\cdot,T)\|_{L^2(\Omega,\mathcal F_T,\mathbb P;L^2(I))}^2.
\end{align}
For the term of $R_4$, by $({\rm d}V)^2=\Theta_1^2|F_2|^2{\rm d}t$ and
\begin{align*}
({\rm d}V_x)^2=&(s\Phi_{1,x}\Theta_1 {\rm d}v^\varepsilon+\Theta_1 {\rm d}{v^\varepsilon_x})^2\nonumber\\
=& \Theta_1^2 |F_{2,x}|^2{\rm d}t+2s\Phi_{1,x}\Theta_1^2 F_2 F_{2,x}{\rm d}t+s^2\Phi_{1,x}^2\Theta_1^2|F_2|^2{\rm d}t,
\end{align*}
we have
\begin{align*}&-\mathbb E\int_{Q_T} R_4{\rm d}x\nonumber\\
=&\frac{1}{2}\mathbb E\int_{Q_T}(x+\varepsilon)^\alpha\big({\rm d}V_x\big)^2{\rm d}x-\frac{1}{2}\mathbb E\int_{Q_T}(x+\varepsilon)^\alpha L_{1,x}^2({\rm d}V)^2{\rm d}x-\mathbb E\int_{Q_T}\frac{1}{2}L_{1,t}({\rm d}V)^2{\rm d}x\nonumber\\
=&\frac{1}{2}\mathbb E\int_{Q_T}(x+\varepsilon)^\alpha\Theta_1^2|F_{2,x}|^2{\rm d}x{\rm d}t+\mathbb E\int_{Q_T}\beta s\lambda(x+\varepsilon)^{\alpha+\beta-1}\Phi_1\Theta_1^2 F_2 F_{2,x}{\rm d}x{\rm d}t\nonumber\\
&-\mathbb E\int_{Q_T}s\lambda(\lambda-t)\Phi_1\Theta_1^2|F_2|^2{\rm d}x{\rm d}t\nonumber\\
\leq & \frac{1}{2}\mathbb E\int_{Q_T}(x+\varepsilon)^\alpha\Theta_1^2|F_{2,x}|^2{\rm d}x{\rm d}t+C(\epsilon)\mathbb E\int_{Q_T}s\Phi_1\Theta_1^2|F_{2,x}|^2{\rm d}x{\rm d}t\nonumber\\
&-\mathbb E\int_{Q_T}s\lambda\left[(1-\epsilon)\lambda-t\right]\Phi_1\Theta_1^2|F_2|^2{\rm d}x{\rm d}t,
\end{align*}
which implies
\begin{align}\label{3.53}
-\mathbb E\int_{Q_T} R_4{\rm d}x\leq C\mathbb E\int_{Q_T}s\Phi_1\Theta_1^2|F_{2,x}|^2{\rm d}x{\rm d}t-\mathbb E\int_{Q_T}s\lambda^2\Phi_1\Theta_1^2|F_2|^2{\rm d}x{\rm d}t.
\end{align}
for sufficiently small $\epsilon$ and sufficiently large $\lambda$ and $s$.

Then substituting (\ref{3.51})-(\ref{3.53}) into (\ref{3.50}) yields that
\begin{align*}
&\mathbb E\int_{Q_T}s\lambda^2\Phi_1\Theta_1^2|F_2|^2{\rm d}x{\rm d}t+\mathbb E\int_{Q_T} s^3\lambda^3(x+\varepsilon)^{2\alpha+3\beta-4}\Phi_1^3\Theta_1^2|v^\varepsilon|^2{\rm d}x {\rm d}t\nonumber\\
&+\mathbb E\int_{Q_T} s \lambda(x+\varepsilon)^{2\alpha+\beta-2}\Phi_1\Theta_1^2|v^\varepsilon_x|^2{\rm d}x{\rm d}t
\end{align*}
\begin{align}\label{3.54}
\leq &C\mathbb E\int_{Q_T}\Theta_1^2 |f_2|^2{\rm d}x{\rm d}t+C\mathbb E\int_{Q_T}s\Phi_1\Theta_1^2|F_{2,x}|^2{\rm d}x{\rm d}t+C(\lambda)\mathbb E\int_0^T s\left[\Theta_1^2|v^\varepsilon_x|^2\right]_{x=1} {\rm d}t\nonumber\\
&+C(\lambda)s^2 e^{C(\lambda)s}\|v^\varepsilon(\cdot,T)\|_{L^2(\Omega,\mathcal F_T,\mathbb P;L^2(I))}^2\nonumber\\
&+C(\lambda)\mathbb E\int_0^T s^2\left[(x+\varepsilon)^{\alpha+\beta-1}\Theta_1^2 (|v^\varepsilon|^2+|F_2|^2)\right]_{x=0}{\rm d}t.
\end{align}

In order to deal with Carleman estimate for nondegenerate part, we use $\Phi_2$ as weight function in Carleman estimate. Letting $L_2=s\Phi_2$, $\Theta_2=e^{L_2}$ and  repeating the above process,   we have the following estimate:
\begin{align}\label{3.55}
&\mathbb E\int_{Q_T}s\lambda^2\Phi_2\Theta_2^2|F_2|^2{\rm d}x{\rm d}t+\mathbb E\int_{Q_T} s^3\lambda^4(x+\varepsilon)^{2\alpha}\eta_{2,x}^4\Phi_2^3\Theta_2^2|v^\varepsilon|^2{\rm d}x {\rm d}t\nonumber\\
&+\mathbb E\int_{Q_T} s\lambda^2(x+\varepsilon)^{2\alpha}\eta_{2,x}^2\Phi_2\Theta_2^2|v^\varepsilon_x|^2{\rm d}x{\rm d}t\nonumber\\
\leq &C\mathbb E\int_{Q_T} s^3\lambda^3(x+\varepsilon)^{\gamma}\Phi_2^3\Theta_2^2|v^\varepsilon|^2{\rm d}x {\rm d}t+C\mathbb E\int_{Q_T} s\lambda(x+\varepsilon)^{\gamma}\Phi_2\Theta_2^2|v^\varepsilon_x|^2{\rm d}x{\rm d}t\nonumber\\
&+C\mathbb E\int_{Q_T}\Theta_2^2 |f_2|^2{\rm d}x{\rm d}t+C\mathbb E\int_{Q_T}s\Phi_2\Theta_2^2|F_{2,x}|^2{\rm d}x{\rm d}t+C(\lambda)s^2 e^{C(\lambda)s}\|v^\varepsilon(\cdot,T)\|_{L^2(\Omega,\mathcal F_T,\mathbb P;L^2(I))}^2\nonumber\\
&+C(\lambda)\mathbb E\int_0^T s^2\left[(x+\varepsilon)^{\alpha}\Theta_2^2 (|v^\varepsilon|^2+|F_2|^2)\right]_{x=0}{\rm d}t
\end{align}
for sufficient large $\lambda$ and $s$ such that $s\geq C(\lambda)$, where $\gamma=\min\{2\alpha-2,\alpha-1\}$.

Now we apply (\ref{3.54}) and (\ref{3.55}) to obtain two Carleman estimates for degenerate part and nondegenerate part, respectively. Obviously, $\tilde v=\chi v^\varepsilon$ satisfies
\begin{align}\label{3.56}\left\{
\begin{array}{ll}
{\rm d}\tilde v^\varepsilon+\left((x+\varepsilon)^\alpha \tilde v^\varepsilon_x\right)_x{\rm d}t=\tilde f_2{\rm d}t+\tilde F_2{\rm d}B(t),&(x,t)\in Q_T,\\
 \tilde v^\varepsilon(1,t)=0,& t\in (0,T),\\
 {\rm and}\ \left\{\begin{array}{ll}
 \tilde v^\varepsilon(0,t)=0 & {\rm for}\ \alpha\in (0,1),\\
 ((x+\varepsilon)^\alpha \tilde v^\varepsilon_x)(0,t)=0 & {\rm for}\ \alpha\in [1,2),
 \end{array}\right. &t\in (0,T).
 \end{array}\right.
\end{align}
where
\begin{align*}
\tilde f_2=-\left((x+\varepsilon)^{\alpha}\chi_x v^\varepsilon\right)_x-(x+\varepsilon)^{\alpha}\chi_x v^\varepsilon_x+\chi f_2,\quad \tilde F_2=\chi F_2.
\end{align*}
Therefore, by applying (\ref{3.54}) to $\tilde v$ and using (\ref{3.42}), we have
\begin{align}\label{3.57}
&\mathbb E\int_{Q_T}s\lambda^2\Phi_1\chi^2|F_2|^2e^{2s\Phi_1}{\rm d}x{\rm d}t+\mathbb E\int_{0}^T\int_0^{x^{(2)}_1} s^3\lambda^3(x+\varepsilon)^{2\alpha+3\beta-4}\Phi^3|v^\varepsilon|^2e^{2s\Phi}{\rm d}x {\rm d}t\nonumber\\
&+\mathbb E\int_{0}^T\int_0^{x^{(2)}_1} s \lambda(x+\varepsilon)^{2\alpha+\beta-2}\Phi|v^\varepsilon_x|^2e^{2s\Phi}{\rm d}x{\rm d}t\nonumber\\
\leq &C\mathbb E\int_{0}^T\int_0^{x^{(2)}_2}|f_2|^2e^{2s\Phi}{\rm d}x{\rm d}t+C\mathbb E\int_{0}^T\int_0^{x^{(2)}_2}s\Phi(|F_{2,x}|^2+|F_2|^2)e^{2s\Phi}{\rm d}x{\rm d}t\nonumber\\
&+C\mathbb E\int_{\omega^{(2)}_T}(|v^\varepsilon|^2+|v^\varepsilon_x|^2)e^{2s\Phi}{\rm d}x{\rm d}t+C(\lambda)s^2 e^{C(\lambda)s}\|v^\varepsilon(\cdot,T)\|_{L^2(\Omega,\mathcal F_T,\mathbb P; L^2(I))}^2\nonumber\\
&+C(\lambda)\mathbb E\int_0^T s^2\left[(x+\varepsilon)^{\alpha+\beta-1} (|v^\varepsilon|^2+|F_2|^2)e^{2s\Phi}\right]_{x=0}{\rm d}t.
\end{align}
Similarly, applying (\ref{3.55}) to $\bar v=(1-\chi) v^\varepsilon$ yields that
\begin{align}
&\mathbb E\int_{Q_T}s\lambda^2\Phi_2(1-\chi)^2|F_2|^2e^{2s\Phi_2}{\rm d}x{\rm d}t+\mathbb E\int_0^T\int_{I\setminus\omega^{(1)}} s^3\lambda^4(x+\varepsilon)^{2\alpha}\Phi_2^3(1-\chi)^2|v^\varepsilon|^2e^{2s\Phi_2}{\rm d}x {\rm d}t\nonumber\\
&+\mathbb E\int_0^T\int_{I\setminus\omega^{(1)}} s\lambda^2(x+\varepsilon)^{2\alpha}\Phi_2(1-\chi)^2|v^\varepsilon_x|^2e^{2s\Phi_2}{\rm d}x{\rm d}t\nonumber\\
\leq &C\mathbb E\int_{Q_T} s^3\lambda^3(x+\varepsilon)^{\gamma}\Phi_2^3(1-\chi)^2|v^\varepsilon|^2e^{2s\Phi_2}{\rm d}x {\rm d}t\nonumber\\
&+C\mathbb E\int_{Q_T} s\lambda(x+\varepsilon)^{\gamma}\Phi_2\big[(1-\chi)^2|v^\varepsilon_x|^2+\chi_x^2|v^\varepsilon|^2\big]e^{2s\Phi_2}{\rm d}x{\rm d}t+C\mathbb E\int_{Q_T} (1-\chi)^2|f_2|^2e^{2s\Phi_2}{\rm d}x{\rm d}t\nonumber\\
&+C\mathbb E\int_{\omega^{(2)}_T}(|v^\varepsilon|^2+|v^\varepsilon_x|^2)e^{2s\Phi_2}{\rm d}x{\rm d}t+C\mathbb E\int_{Q_T}s\Phi_2\big[(1-\chi)^2|F_{2,x}|^2+\chi_x^2|F_{2}|^2\big]e^{2s\Phi_2}{\rm d}x{\rm d}t\nonumber\\
&+C(\lambda)s^2 e^{C(\lambda)s}\|v^\varepsilon(\cdot,T)\|_{L^2(\Omega,\mathcal F_T,\mathbb P;L^2(I))}^2.
\end{align}
Obviously, $(1-\chi)^2(x+\varepsilon)^{\gamma}\leq C(1-\chi)^2(x+\varepsilon)^{2\alpha}$ in $Q_T$, where $C$ is not depending on $\varepsilon$. Then, we further obtain for any $\epsilon>0$  that
\begin{align}\label{3.59}
&\mathbb E\int_{Q_T}s\lambda^2\Phi_2(1-\chi)^2|F_2|^2e^{2s\Phi_2}{\rm d}x{\rm d}t+\mathbb E\int_0^T\int_{x_2^{(2)}}^1 s^3\lambda^4\Phi^3(x+\varepsilon)^{2\alpha+3\beta-4} |v^\varepsilon|^2e^{2s\Phi}{\rm d}x {\rm d}t\nonumber\\
&+\epsilon\mathbb E\int_0^T\int_{x_2^{(2)}}^1 s\lambda^2\Phi(x+\varepsilon)^{2\alpha+\beta-2}|v^\varepsilon_x|^2e^{2s\Phi}{\rm d}x{\rm d}t\nonumber\\
\leq &C\mathbb E\int_{0}^T\int_{x_1^{(2)}}^1 |f_2|^2e^{2s\Phi}{\rm d}x{\rm d}t+C\mathbb E\int_{0}^T\int_{x_1^{(2)}}^1 s\Phi(|F_{2,x}|^2+|F_2|^2)e^{2s\Phi}{\rm d}x{\rm d}t\nonumber\\
&+C(\lambda)s^2 e^{C(\lambda)s}\|v^\varepsilon(\cdot, T)\|_{L^2(\Omega,\mathcal F_T,\mathbb P;L^2(I))}^2\nonumber\\
&+C\mathbb E\int_{\omega^{(1)}_T}(s^3\lambda^4\Phi^3|v^\varepsilon|^2+\epsilon s\lambda^2\Phi|v^\varepsilon_x|^2)e^{2s\Phi}{\rm d}x{\rm d}t
\end{align}
for sufficient large $s$ and $\lambda$.

Combining  (\ref{3.57}) and (\ref{3.59}) and adding to both sides of
the inequality the term
\begin{align*}
\mathbb E\int_0^T\int_{x_1^{(2)}}^{x_2^{(2)}}s^3\lambda^4\Phi^3(x+\varepsilon)^{2\alpha+3\beta-4} |v^\varepsilon|^2e^{2s\Phi}{\rm d}x{\rm d}t+\epsilon\mathbb E\int_0^T\int_{x_1^{(2)}}^{x_2^{(2)}}s\lambda^2\Phi(x+\varepsilon)^{2\alpha+\beta-2}
|v^\varepsilon_x|^2e^{2s\Phi}{\rm d}x{\rm d}t,
\end{align*}
 we obtain
\begin{align}
&\mathbb E\int_{Q_T}s\lambda^2\left[\Phi_1\chi^2e^{2s\Phi_1}+\Phi_2(1-\chi)^2e^{2s\Phi_2}\right]|F_2|^2{\rm d}x{\rm d}t\nonumber\\
&+\mathbb E\int_{Q_T} s^3\lambda^4\Phi^3(x+\varepsilon)^{2\alpha+3\beta-4} |v^\varepsilon|^2e^{2s\Phi}{\rm d}x {\rm d}t+\epsilon\mathbb E\int_{Q_T} s\lambda^2\Phi(x+\varepsilon)^{2\alpha+\beta-2}|v^\varepsilon_x|^2e^{2s\Phi}{\rm d}x{\rm d}t\nonumber\\
\leq &C\mathbb E\int_{Q_T}\lambda |f_2|^2e^{2s\Phi}{\rm d}x{\rm d}t+C\mathbb E\int_{Q_T} s\lambda\Phi\big(|F_{2,x}|^2+|F_{2}|^2\big)e^{2s\Phi}{\rm d}x{\rm d}t\nonumber\\
&+C\mathbb E\int_{\omega^{(1)}_T}(s^3\lambda^4\Phi^3|v^\varepsilon|^2+\epsilon s\lambda^2\Phi|v^\varepsilon_x|^2)e^{2s\Phi}{\rm d}x{\rm d}t+C(\lambda)s^2 e^{C(\lambda)s}\|v^\varepsilon(\cdot,T)\|_{L^2(\Omega,\mathcal F_T,\mathbb P;L^2(I))}^2\nonumber\\
&+C(\lambda)\mathbb E\int_0^T s^2\left[(x+\varepsilon)^{\alpha+\beta-1} (|v^\varepsilon|^2+|F_2|^2)e^{2s\Phi}\right]_{x=0}{\rm d}t.
\end{align}
Noticing that $\Phi_1\chi^2e^{2s\Phi_1}+\Phi_2(1-\chi)^2e^{2s\Phi_2}>C\Phi e^{2s\Phi}$ in $Q_T$,  we further have
\begin{align}\label{3.61}
&\mathbb E\int_{Q_T}s\lambda^2\Phi|F_2|^2e^{2s\Phi}{\rm d}x{\rm d}t+\mathbb E\int_{Q_T} s^3\lambda^4\Phi^3(x+\varepsilon)^{2\alpha+3\beta-4} |v^\varepsilon|^2e^{2s\Phi}{\rm d}x {\rm d}t\nonumber\\
&+\epsilon\mathbb E\int_{Q_T} s\lambda^2\Phi(x+\varepsilon)^{2\alpha+\beta-2}|v^\varepsilon_x|^2e^{2s\Phi}{\rm d}x{\rm d}t\nonumber\\
\leq &C\mathbb E\int_{Q_T} \lambda |f_2|^2e^{2s\Phi}{\rm d}x{\rm d}t+C\mathbb E\int_{Q_T} s\lambda\Phi|F_{2,x}|^2e^{2s\Phi}{\rm d}x{\rm d}t\nonumber\\
&+C\mathbb E\int_{\omega^{(1)}_T}(s^3\lambda^4\Phi^3|v^\varepsilon|^2+\epsilon s\lambda^2\Phi|v^\varepsilon_x|^2)e^{2s\Phi}{\rm d}x{\rm d}t+C(\lambda)s^2 e^{C(\lambda)s}\|v^\varepsilon(\cdot,T)\|_{L^2(\Omega,\mathcal F_T,\mathbb P;L^2(I))}^2\nonumber\\
&+C(\lambda)\mathbb E\int_0^T s^2\left[(x+\varepsilon)^{\alpha+\beta-1} (|v^\varepsilon|^2+|F_2|^2)e^{2s\Phi}\right]_{x=0}{\rm d}t.
\end{align}

Similar to Lemma 3.4,  we have the following Cacciopoli inequality for forward stochastic degenerate parabolic equation:
\begin{align}\label{3.62}
&\mathbb E \int_{\omega^{(1)}_T} \Phi|v^\varepsilon_x|^2e^{2s\Phi}{\rm d}x{\rm d}t\nonumber\\
\leq& C(\lambda)\mathbb E \int_{\omega_T}s^2\Phi^3|v^\varepsilon|^2e^{2s\Phi}{\rm d}x{\rm d}t+C\mathbb E\int_{Q_T} s^{-1}\Phi|f_2|^2e^{2s\Phi}{\rm d}x{\rm d}t\nonumber\\
&+C\mathbb E \int_{Q_T} \Phi|F_2|^2e^{2s\Phi}{\rm d}x{\rm d}t.
\end{align}
Finally, substituting (\ref{3.62}) into (\ref{3.61}) and choosing $\epsilon$ sufficiently small, we can absorb the term of $F_2$ on the right-hand side of (\ref{3.62}) and then obtain (\ref{3.43}). This completes the proof of Theorem 3.5. \hfill $\Box$

\section{Null controllability}
\setcounter{equation}{0}
In this section, we will apply Theorem 3.1 to prove the null controllability result for the forward stochastic degenerate parabolic equation (\ref{1.3}), i.e. the following Theorem 4.1.
\begin{thm}
Let $\alpha\in (0,\frac{1}{2})$ and $a,b,c\in L^\infty_\mathcal F(0,T;L^\infty(I))$. Then for any $y_0\in L^2(\Omega,\mathcal F_0,\mathbb P;$ $L^2(G))$, there exists a pair of controls $(g,G)\in L^2_\mathcal F(0,T;L^2(I))\times L^2_\mathcal F(0,T;L^2(I))$ such that the solution $y$ of (\ref{1.3}) satisfies $y(x,T)=0$ in $I$, $\mathbb P$-a.s.
\end{thm}

Since the system (\ref{1.3}) is degenerate, we first transfer to study a uniform null
controllability in $\varepsilon$ for a nondegenerate approximate system. More precisely, letting $0<\varepsilon<1$, we consider
\begin{align}\label{4.1}\hspace{-0.1cm}\left\{
\begin{array}{ll}
{\rm d}y^\varepsilon-\left((x+\varepsilon)^\alpha y^\varepsilon_x\right)_x{\rm d}t=(a y^\varepsilon_{x}+b y^\varepsilon+g^\varepsilon{\bf 1}_\omega){\rm d}t+(cy^\varepsilon+G^\varepsilon){\rm d}B(t),&(x,t)\in Q_T,\\
 y^\varepsilon(0,t)=y^\varepsilon(1,t)=0,& t\in (0,T),\\
 y^\varepsilon(x,0)=y_0^\varepsilon(x),&x\in I,
 \end{array}\right.
\end{align}
where
\begin{align}
y^\varepsilon_0\rightarrow y_0\quad {\rm in}\ L^2(\Omega,\mathcal F_0,\mathbb P;L^2(I)).
\end{align}
  It is well known that the key ingredient for  studying the null controllability is to obtain observation inequality for the corresponding adjoint equation. An important tool is Carleman estimate, in whose proof  the main difficulty is how to deal with the first order term  in the stochastic degenerate parabolic system.
In order to use the terms on the left-hand side of Carleman estimate to absorb  this term directly, we need $x^{-\frac{\alpha}{2}}a\in L^\infty_{\mathcal F}(0,T;L^\infty(I))$, which means that the coefficient $a$ of the first order
term goes to zero at some polynomial rate as $x\rightarrow 0$. More reasonable condition is $a\in L^\infty_{\mathcal F}(0,T;L^\infty(I))$. For this condition on $a$, we will apply a duality technique to establish a new Carleman estimate for the stochastic degenerate parabolic equation with convection term.

In next subsection we first prove a Carleman estimate for the following corresponding adjoint system of (\ref{4.1}):
\begin{align}\label{4.3}\left\{
\begin{array}{ll}
{\rm d}z+\left((x+\varepsilon)^\alpha z_x\right)_x{\rm d}t=\big((a z)_{x}-b z-cZ\big){\rm d}t+Z{\rm d}B(t),&(x,t)\in Q_T,\\
 z(0,t)=z(1,t)=0,& t\in (0,T),\\
 z(x,T)=z_T(x),&x\in I.
 \end{array}\right.
\end{align}
 Next based on this Carleman estimate, we obtain observation inequality and then prove the null controllability result, i.e. Theorem 4.1.

\subsection{Carleman estimate for a backward stochastic degenerate equation with convection term}

Our main result in this subsection is the following estimate, whose proof is based on a duality argument introduced by  Imanuvilov and Yamamoto [\ref{Imanuvilov2001}] for deterministic parabolic equation, or introduced by Liu [\ref{Liu2014ECOCV}] or Yan [\ref{Yan2018JMAA}] for stochastic parabolic equation.

\begin{thm}
 Let $\alpha\in \left(0,\frac{1}{2}\right)$, $a, b, c\in L^\infty_\mathcal F(0,T; L^\infty(I))$ and $z_T\in L^2(\Omega,\mathcal F_T,\mathbb P; L^2(I))$.
Then for any $\varepsilon\in (0,\nu)$ with a sufficiently small $\nu>0$,
 there exist positive constants $\lambda_3=\lambda_3(\omega,I,T$, $\alpha,M,\nu)$, $s_3=s_3(\omega,I,T$, $\alpha,M,\nu,\lambda)$ and
$C=C(\omega,I,T,\alpha,M,\nu,\lambda)$ such that
\begin{align}\label{4.4}
&\mathbb E\int_{Q_T}s^3\xi^3|z|^2e^{2s\varphi}{\rm d}x{\rm d}t+\mathbb E\int_{Q_T}s\xi(x+\varepsilon)^{\alpha}|z_x|^2e^{2s\varphi}{\rm d}x{\rm d}t\nonumber\\
\leq& C\mathbb E\int_{Q_T}s^{2}\xi^{2}|Z|^2e^{2s\varphi}{\rm d}x{\rm
d}t+C\mathbb E\int_{\omega_T}s^3\xi^3|z|^2e^{2s\varphi}{\rm d}x{\rm
d}t
\end{align}
for all $\lambda\geq \lambda_3$, $s\geq s_3$ and all $z\in \mathcal H^1$ satisfying (\ref{4.3}).

\end{thm}

In order to prove Theorem 4.2, we consider the following controlled forward stochastic parabolic equation:
\begin{align}\label{4.5}\left\{\begin{array}{ll}
{\rm d} w-\left((x+\varepsilon)^\alpha { w}_x\right)_x{\rm d}t=\left(s^3\xi^3 ze^{2s\varphi}+ h{\mathbf 1}_{\omega}\right){\rm d}t+H{\rm d}B(t),&(x,t)\in Q_T,\\
w(0,t)=w(1,t)=0,&t\in (0,T),\\
w(x,0)=0, &x\in I,\end{array}\right.
\end{align}
where $(h,H)\in L^2_\mathcal F(0,T;L^2(\omega))\times L^2_\mathcal F(0,T;L^2(I))$ is a pair control. Then, we have the following controllability
result, whose proof will be put in the Appendix.

\begin{lem}
Let $\alpha\in \big(0,\frac{1}{2}\big)$. Then for any $\varepsilon\in (0,1)$, there exists a pair of controls $(h, H)\in L_\mathcal F^2(0,T;L^2(\omega))\times L_\mathcal F^2(0,T;L^2(I))$ such that (\ref{4.5}) admits a solution $w\in \mathcal H^1$ corresponding to $(h, H)$ satisfying $w(x,T)=0$ in $I$, P-$a.s.$ Moreover, there exists a positive constant $C=C(\omega,I,T,\alpha,M)$ such that
\begin{align}\label{4.6}
&\mathbb E\int_{Q_T} |w|^2e^{-2s\varphi}{\rm d}x{\rm d}t+\mathbb E\int_{Q_T}s^{-2}\xi^{-2}(x+\varepsilon)^\alpha|w_{x}|^2e^{-2s\varphi}{\rm d}x{\rm d}t\nonumber\\
&+\mathbb E\int_{\omega_T} s^{-3}\xi^{-3} |h|^2e^{-2s\varphi}{\rm d}x{\rm d}t+\mathbb E \int_{Q_T}s^{-2}\xi^{-2} |H|^2e^{-2s\varphi}{\rm d}x{\rm d}t\nonumber\\
\leq &C \mathbb E \int_{Q_T}s^3\xi^3 |z|^2e^{2s\varphi}{\rm d}x{\rm d}t.
\end{align}
\end{lem}

\vspace{2mm}

Now we prove Theorem 4.2.

\vspace{2mm}

\noindent{\bf Proof of Theorem 4.2.}\  By Lemma 4.3, we know that  there exists a pair of controls $(h,H)$ such that the solution $w$ of (\ref{4.5}) corresponding to $(h,H)$ satisfies $w(x,T)=0$ in $I$, P-$a.s$. Then by using It\^{o} formula and integrating by parts, we obtain the following duality between $w$ and $z$:
\begin{align}
&\mathbb E\int_{Q_T}\left(s^3\xi^3 |z|^2 e^{2s\varphi}+z h|_{\omega}\right){\rm d}x{\rm d}t\nonumber\\
=&\mathbb E\int_{Q_T}azw_x {\rm d}x{\rm d}t+\mathbb E\int_{Q_T} (bzw+cZw){\rm d}x{\rm d}t-\mathbb E\int_{Q_T} Z H{\rm d}x{\rm d}t.
\end{align}
By Young's inequality, we further find that
\begin{align}\label{4.8}
&\mathbb E\int_{Q_T}s^3\xi^3 |z|^2 e^{2s\varphi}{\rm d}x{\rm d}t\nonumber\\
\leq &\epsilon \mathbb E\int_{\omega_T} s^{-3}\xi^{-3} |h|^2e^{-2s\varphi}{\rm d}x{\rm d}t +\epsilon\mathbb E\int_{Q_T}s^{-2}\xi^{-2}(x+\varepsilon)^\alpha|w_{x}|^2e^{-2s\varphi}{\rm d}x{\rm d}t\nonumber\\
&+\epsilon\mathbb E\int_{Q_T} |w|^2e^{-2s\varphi}{\rm d}x{\rm d}t+\epsilon\mathbb E \int_{Q_T}s^{-2}\xi^{-2} |H|^2e^{-2s\varphi}{\rm d}x{\rm d}t\nonumber\\
&+C(\epsilon)\mathbb E\int_{\omega_T}s^3\xi^3 |z|^2 e^{2s\varphi}{\rm d}x{\rm d}t+C(\epsilon)\mathbb E\int_{Q_T}s^{2}\xi^{2}(x+\varepsilon)^{-\alpha}|z|^2e^{2s\varphi}{\rm d}x{\rm d}t\nonumber\\
&+C(\epsilon)\mathbb E\int_{Q_T} |z|^2e^{2s\varphi}{\rm d}x{\rm d}t+C(\epsilon)\mathbb E \int_{Q_T}s^{2}\xi^{2} |Z|^2e^{2s\varphi}{\rm d}x{\rm d}t.
\end{align}
Substituting (\ref{4.6}) into (\ref{4.8}) and choosing $\epsilon$ sufficiently small, we obtain
 \begin{align}\label{4.9}
\mathbb E\int_{Q_T}s^3\xi^3 |z|^2 e^{2s\varphi}{\rm d}x{\rm d}t\leq & C\mathbb E\int_{Q_T}s^{2}\xi^{2}(x+\varepsilon)^{-\alpha}|z|^2e^{2s\varphi}{\rm d}x{\rm d}t+C\mathbb E \int_{Q_T}s^{2}\xi^{2} |Z|^2e^{2s\varphi}{\rm d}x{\rm d}t\nonumber\\
&+C\mathbb E\int_{\omega_T}s^3\xi^3 |z|^2 e^{2s\varphi}{\rm d}x{\rm d}t.
\end{align}

Now we estimate $\mathbb E\int_{Q_T}s\xi(x+\varepsilon)^{\alpha}
|z_x|^2e^{2s\varphi}{\rm d}x{\rm d}t$. To do this, we use It\^{o} formula again and the equation of $z$ to obtain
\begin{align}\label{}
&2\mathbb E\int_{Q_T}s\xi(x+\varepsilon)^{\alpha}|z_x|^2e^{2s\varphi}{\rm d}x{\rm d}t\nonumber\\
=&-\mathbb E\int_{Q_T}s(\xi e^{2s\varphi})_t|z|^2{\rm d}x{\rm d}t-2\mathbb E\int_{Q_T}s\xi(x+\varepsilon)^\alpha(e^{2s\varphi})_x zz_x{\rm d}x{\rm d}t\nonumber\\
&+2\mathbb E\int_{Q_T}s\xi a z \left(z e^{2s\varphi}\right)_x{\rm d}x{\rm d}t+2\mathbb E\int_{Q_T}s\xi z(bz+cZ)e^{2s\varphi} {\rm d}x{\rm d}t-\mathbb E\int_{Q_T} s\xi|Z|^2 e^{2s\varphi}{\rm d}x{\rm d}t\nonumber\\
\leq & C\mathbb E\int_{Q_T}s^3\xi^3 |z|^2 e^{2s\varphi}{\rm d}x{\rm d}t+\mathbb E\int_{Q_T}s\xi(x+\varepsilon)^{\alpha}|z_x|^2e^{2s\varphi}{\rm d}x{\rm d}t\nonumber\\
&+C\mathbb E\int_{Q_T}s\xi(x+\varepsilon)^{-\alpha}|z|^2e^{2s\varphi}{\rm d}x{\rm d}t+C\mathbb E\int_{Q_T}s\xi|Z|^2e^{2s\varphi}{\rm d}x{\rm d}t,
\end{align}
which implies
\begin{align}\label{4.11}
&\mathbb E\int_{Q_T}s\xi(x+\varepsilon)^{\alpha}|z_x|^2e^{2s\varphi}{\rm d}x{\rm d}t\nonumber\\
\leq & C\mathbb E\int_{Q_T}s^3\xi^3 |z|^2 e^{2s\varphi}{\rm d}x{\rm d}t+C\mathbb E\int_{Q_T}s\xi(x+\varepsilon)^{-\alpha}|z|^2e^{2s\varphi}{\rm d}x{\rm d}t\nonumber\\
&+C\mathbb E\int_{Q_T}s\xi|Z|^2e^{2s\varphi}{\rm d}x{\rm d}t.
\end{align}
 Combining (\ref{4.9}) and (\ref{4.11}) yields that
   \begin{align}\label{4.12}
&\mathbb E\int_{Q_T}s^3\xi^3 |z|^2 e^{2s\varphi}{\rm d}x{\rm d}t+\mathbb E\int_{Q_T}s\xi(x+\varepsilon)^{\alpha}|z_x|^2e^{2s\varphi}{\rm d}x{\rm d}t\nonumber\\
\leq & C\mathbb E\int_{Q_T}s^{2}\xi^{2}(x+\varepsilon)^{-\alpha}|z|^2e^{2s\varphi}{\rm d}x{\rm d}t+C\mathbb E \int_{Q_T}s^{2}\xi^{2} |Z|^2e^{2s\varphi}{\rm d}x{\rm d}t\nonumber\\
&+C\mathbb E\int_{\omega_T}s^3\xi^3 |z|^2 e^{2s\varphi}{\rm d}x{\rm d}t.
\end{align}
Applying Young's equality  and Hardy-Poincar\'{e} inequality (\ref{3.5}), we obtain
\begin{align}\label{4.13}
&\mathbb E\int_{Q_T}s^{2}\xi^{2}(x+\varepsilon)^{-\alpha}|z|^2e^{2s\varphi}{\rm d}x{\rm d}t\nonumber\\
\leq &\epsilon\mathbb E\int_{Q_T}s^{3}\xi^{3}|z|^2e^{2s\varphi}{\rm d}x{\rm d}t+C(\epsilon)\mathbb E\int_{Q_T}s\xi (x+\varepsilon)^{-2\alpha}|z|^2e^{2s\varphi}{\rm d}x{\rm d}t\nonumber\\
\leq &\epsilon\mathbb E\int_{Q_T}s^{3}\xi^{3}|z|^2e^{2s\varphi}{\rm d}x{\rm d}t+C(\varepsilon)\mathbb E\int_{Q_T} (x+\varepsilon)^{2-2\alpha}\left(s\xi|z_x|^2+s^3\xi^3|z|^2\right)e^{2s\varphi}{\rm d}x{\rm d}t.
\end{align}
Choosing $\epsilon$ sufficiently small and substituting (\ref{4.13}) into (\ref{4.12}), we find that
\begin{align}\label{4.14}
&\mathbb E\int_{Q_T}s^3\xi^3 |z|^2 e^{2s\varphi}{\rm d}x{\rm d}t+\mathbb E\int_{Q_T}s\xi(x+\varepsilon)^{\alpha}|z_x|^2e^{2s\varphi}{\rm d}x{\rm d}t\nonumber\\
\leq & C\mathbb E \int_{Q_T}s^{2}\xi^{2} |Z|^2e^{2s\varphi}{\rm d}x{\rm d}t+C\mathbb E\int_{\omega_T}s^3\xi^3 |z|^2 e^{2s\varphi}{\rm d}x{\rm d}t\nonumber\\
&+C\mathbb E\int_{Q_T} (x+\varepsilon)^{2-2\alpha}\left(s\xi|z_x|^2+s^3\xi^3|z|^2\right)e^{2s\varphi}{\rm d}x{\rm d}t.
\end{align}

The remainder of the proof is to eliminate the last term on the right-hand side of (\ref{4.14}).  In order to overcome the degeneracy in this term, we transfer to consider the equation of $z$ in a interval outside of $x=0$. For some given $0<\nu<\frac{x^{(1)}_1}{3}$, we set $I_\nu=(\nu,1)$ and $Q_{\nu,T}=I_\nu\times (0,T)$. Further we introduce a cut-off function $\rho\in C^2(\bar{I})$ such that
$0\leq \rho(x)\leq 1$ for $x\in  I$, $\rho(x)\equiv1$ for $x\in
(3\nu,1)$ and $\rho(x)\equiv0$ for $x\in (0,2\nu)$. Additionally, we choose a weight function $\tilde \varphi$ such that
\begin{align*}
\tilde\varphi(x,t)=\frac{e^{\lambda\tilde\eta(x)}-e^{2\lambda M}}{t^2(T-t)^2},\quad (x,t)\in Q_{\nu,T},
\end{align*}
where $\tilde \eta\in C^2(\overline I_{\nu})$ satisfies $\tilde \eta>0$ in $I_{\nu}$ and
\begin{align}\tilde\eta(x)=\left\{
\begin{array}{ll}
0,&x=\nu,\\
{\rm smooth},&x\in (\nu,2\nu),\\
\eta_1(x), &x\in (2\nu,x^{(2)}_1),\\
\eta_1(x)=\eta_2(x),&x\in (x^{(2)}_1,x^{(2)}_2),\\
\eta_2(x), &x\in (x^{(2)}_2,1).
\end{array}\right.
\end{align}
Then we easily see that
\begin{align*}
\tilde \eta>0,\quad x\in I_\nu,\quad \tilde \eta(\nu)=\tilde \eta(1)=0\quad{\rm
and}\quad |\tilde \eta_{x}(x)|>0,\quad  x\in
\overline{I_\nu\setminus\omega^{(1)}}
\end{align*}
and
\begin{align}\label{4.16}
\tilde \varphi(x,t)=\varphi(x,t),\quad (x,t)\in (2\nu,1)\times (0,T).
\end{align}
 Letting $\tilde z =\rho z$, we now consider
\begin{align}\left\{
\begin{array}{ll}
{\rm d}\tilde z+\left((x+\varepsilon)^\alpha \tilde z_x\right)_x{\rm d}t=\big((a \tilde z)_{x}-b \tilde z-c\tilde Z+\tilde f\big){\rm d}t+\tilde Z{\rm d}B(t),&(x,t)\in Q_{\nu,T},\\
 z(\nu,t)=z(1,t)=0,& t\in (0,T),\\
 z(x,T)=z_T(x),&x\in I_{\nu},
 \end{array}\right.
\end{align}
 where
\begin{align*}
\tilde f=\big((x+\varepsilon)^\alpha\rho_x z\big)_x+(x+\varepsilon)^\alpha \rho_x z_x-a\rho_x z,\quad \tilde Z=\rho Z.
\end{align*}
By the Carleman estimate for  stochastic nondegenerate parabolic equation, e.g. Theorem 6.1 in  [\ref{TangSIAM2009}], we obtain that there exists a constant $C$ depending on $\omega, I, T, \alpha$ and $\nu$, but independent of $\varepsilon$ such that
\begin{align}
&\mathbb E\int_{Q_{\nu,T}}s\xi |\tilde z_x|^2e^{2s\tilde\varphi}{\rm d}x{\rm d}t+\mathbb E\int_{Q_{\nu,T}}s^3\xi^3 |\tilde z|^2e^{2s\tilde\varphi}{\rm d}x{\rm d}t\nonumber\\
\leq & C(\nu) \mathbb E \int_{Q_{\nu,T}}(|\tilde Z|^2+|\tilde f|^2)e^{2s\tilde\varphi}{\rm d}x{\rm d}t+C(\nu)\mathbb E \int_{Q_{\nu,T}}s^2\xi^2|\tilde Z|^2e^{2s\tilde\varphi}{\rm d}x{\rm d}t\nonumber\\
&+C(\nu)\mathbb E \int_{\omega^{(1)}_{T}}s^3\xi^3|\tilde z|^2e^{2s\tilde\varphi}{\rm d}x{\rm d}t\nonumber\\
\leq &C(\nu) \mathbb E \int_{Q_{\nu,T}}(|\rho_x|^2+|\rho_{xx}|^2)(| z|^2+| z_x|^2)e^{2s\tilde\varphi}{\rm d}x{\rm d}t+C(\nu)\mathbb E \int_{Q_{\nu,T}}s^2\xi^2|\tilde Z|^2e^{2s\tilde\varphi}{\rm d}x{\rm d}t\nonumber\\
&+C(\nu)\mathbb E \int_{\omega^{(1)}_{T}}s^3\xi^3|\tilde z|^2e^{2s\tilde\varphi}{\rm d}x{\rm d}t.
\end{align}
Using the definition of $\rho$ and (\ref{4.16}), we further obtain
\begin{align}\label{4.19}
&\mathbb E\int_0^T\int_{3\nu}^1s\xi |z_x|^2e^{2s\varphi}{\rm d}x{\rm d}t+\mathbb E\int_0^T\int_{3\nu}^1 s^3\xi^3 |z|^2e^{2s\varphi}{\rm d}x{\rm d}t\nonumber\\
\leq &C(\nu) \mathbb E \int_0^T\int_{2\nu}^{3\nu}(|z|^2+| z_x|^2)e^{2s\varphi}{\rm d}x{\rm d}t+C(\nu)\mathbb E \int_0^T\int_{2\nu}^1s^2\xi^2| Z|^2e^{2s\varphi}{\rm d}x{\rm d}t\nonumber\\
&+C(\nu)\mathbb E \int_{\omega^{(1)}_{T}}s^3\xi^3|z|^2e^{2s \varphi}{\rm d}x{\rm d}t.
\end{align}
On the other hand,  we easily obtain that
\begin{align}\label{4.20}
&\mathbb E\int_{Q_T} (x+\varepsilon)^{2-2\alpha}\left(s\xi|z_x|^2+s^3\xi^3|z|^2\right)e^{2s\varphi}{\rm d}x{\rm d}t\nonumber\\
\leq & \mathbb E\int_{0}^T\int_{0}^{3\nu} (x+\varepsilon)^{2-2\alpha}\left(s\xi|z_x|^2+s^3\xi^3|z|^2\right)e^{2s\varphi}{\rm d}x{\rm d}t\nonumber\\
&+\mathbb E\int_{0}^T\int_{3\nu}^{1} (x+\varepsilon)^{2-2\alpha}\left(s\xi|z_x|^2+s^3\xi^3|z|^2\right)e^{2s\varphi}{\rm d}x{\rm d}t\nonumber\\
\leq& (4\nu)^{2-3\alpha}\mathbb E\int_{0}^T\int_{0}^{3\nu}\left(s\xi(x+\varepsilon)^\alpha|z_x|^2+s^3\xi^3|z|^2\right)e^{2s\varphi}{\rm d}x{\rm d}t\nonumber\\
&+2^{2-2\alpha}\mathbb E\int_{0}^T\int_{3\nu}^{1} \left(s\xi|z_x|^2+s^3\xi^3|z|^2\right)e^{2s\varphi}{\rm d}x{\rm d}t
\end{align}
for $\nu\in \left(0,\frac{1}{4}\right)$ and $\varepsilon\in (0,\nu)$. Then from (\ref{4.19}) and (\ref{4.20}) it follows that
\begin{align}\label{4.21}
&\mathbb E\int_{Q_T}(x+\varepsilon)^{2-2\alpha}\left(s\xi|z_x|^2+s^3\xi^3|z|^2\right)e^{2s\varphi}{\rm d}x{\rm d}t\nonumber\\
\leq& (4\nu)^{2-3\alpha}\mathbb E\int_{0}^T\int_{0}^{3\nu}\left(s\xi(x+\varepsilon)^\alpha|z_x|^2+s^3\xi^3|z|^2\right)e^{2s\varphi}{\rm d}x{\rm d}t\nonumber\\
&+C(\nu) \mathbb E \int_0^T\int_{2\nu}^{3\nu}(|z|^2+| z_x|^2)e^{2s\varphi}{\rm d}x{\rm d}t+C(\nu)\mathbb E \int_{Q_T}s^2\xi^2| Z|^2e^{2s\varphi}{\rm d}x{\rm d}t\nonumber\\
&+C(\nu)\mathbb E \int_{\omega^{(1)}_{T}}s^3\xi^3|z|^2e^{2s \varphi}{\rm d}x{\rm d}t.
\end{align}
By substituting (\ref{4.21}) into (\ref{4.14}), we obtain
\begin{align}\label{4.22}
&\mathbb E\int_{Q_T}s^3\xi^3 |z|^2 e^{2s\varphi}{\rm d}x{\rm d}t+\mathbb E\int_{Q_T}s\xi(x+\varepsilon)^{\alpha}|z_x|^2e^{2s\varphi}{\rm d}x{\rm d}t\nonumber\\
\leq & C(\nu)\mathbb E \int_{Q_T}s^{2}\xi^{2} |Z|^2e^{2s\varphi}{\rm d}x{\rm d}t+C(\nu)\mathbb E\int_{\omega_T}s^3\xi^3 |z|^2 e^{2s\varphi}{\rm d}x{\rm d}t\nonumber\\
&+(4\nu)^{2-3\alpha}C\mathbb E\int_{Q_T}\left(s\xi(x+\varepsilon)^\alpha|z_x|^2+s^3\xi^3|z|^2\right)e^{2s\varphi}{\rm d}x{\rm d}t\nonumber\\
&+C(\nu) \mathbb E \int_{Q_T}|z|^2e^{2s\varphi}{\rm d}x{\rm d}t+\frac{C(\nu)}{(2\nu)^\alpha} \mathbb E \int_{Q_T}(x+\varepsilon)^\alpha| z_x|^2e^{2s\varphi}{\rm d}x{\rm d}t.
\end{align}
Finally, choosing $\nu$ sufficiently small such that $(4\nu)^{2-3\alpha}C\leq \frac{1}{4}$ and then $s$ sufficiently large such that $\frac{1}{2}s\min_{t\in [0,T]}\xi>\max\big\{C(\nu),\frac{C(\nu)}{(2\nu)^\alpha}\big\}$, we can absorb the last three terms on the right-hand side of (\ref{4.22}) and then obtain (\ref{4.4}). This completes the proof of Theorem 4.2. \hfill$\Box$

\subsection{Proof of Theorem 4.1}

In this section, we will show the null controllability result for system
(\ref{1.3}), i.e. Theorem 4.1. To do this, we first prove the following observation inequality.
\begin{lem}
Let  $\alpha\in \big(0,\frac{1}{2}\big)$, $a, b,c\in L^\infty_\mathcal F(0,T; L^\infty(I))$ and $z_T\in L^2(\Omega,\mathcal F_T,\mathbb P; L^2(I))$. Then for any $\varepsilon\in (0,\nu)$ with a sufficiently small $\nu>0$, there exist positive constant $C=C(\omega,I,T,\alpha,M,\nu)$ such that the solution $z$ of the adjoint
system (\ref{4.3})  satisfies
\begin{align}\label{4.23}
\mathbb E\int_{I} |z|^2(x,0){\rm
d}x\leq C \int_{Q_T} |Z|^2{\rm
d}x+C\int_{\omega_T}|z|^2{\rm d}x{\rm d}t.
\end{align}
\end{lem}

{\noindent\bf Proof.}\ By It\^{o} formula, we obtain for
$0\leq\tau<\tilde\tau\leq T$ that
\begin{align}\label{4.24}
&\mathbb E\int_{I}|z|^2(x,\tau){\rm d}x+2\mathbb E\int_{\tau}^{\tilde \tau}\int_{I} (x+\varepsilon)^{\alpha}|z_{x}|^2{\rm d}x{\rm d}t+\mathbb E\int_{\tau}^{\tilde
\tau}\int_{I} |Z|^2{\rm d}x{\rm d}t
\nonumber\\
= &\mathbb E\int_{I}|z|^2(x,\tilde
\tau){\rm d}x+2\mathbb E\int_{\tau}^{\tilde
\tau}\int_{I}az_xz{\rm d}x{\rm d}t+2\mathbb E\int_{\tau}^{\tilde
\tau}\int_{I} b|z|^2{\rm d}x{\rm d}t\nonumber\\
&+2\mathbb E\int_{\tau}^{\tilde
\tau}\int_{I} czZ{\rm d}x{\rm d}t.
\end{align}
Similar to (\ref{2.14}), we have
\begin{align}\label{4.25}
\mathbb E\int_{\tau}^{\tilde
\tau}\int_{I}az_xz{\rm d}x{\rm d}t\leq& \frac{1}{2} \mathbb E\int_{\tau}^{\tilde
\tau}\int_{I}(x+\varepsilon)^\alpha |z_x|^2{\rm d}x{\rm d}t+C\mathbb E\int_{\tau}^{\tilde
\tau}\int_{I}|z|^2{\rm d}x{\rm d}t.
\end{align}
Substituting (\ref{4.25}) into (\ref{4.24}) yields that
\begin{align}
&\mathbb E\int_{I}|z|^2(x,\tau){\rm d}x\leq \mathbb E\int_{I}|z|^2(x,\tilde
\tau){\rm d}x+C\mathbb E\int_\tau^{\tilde\tau}\int_I |z|^2{\rm d}x{\rm d}t,\quad 0\leq \tau<\tilde\tau\leq T.
\end{align}
Then applying Gronwall inequality yields that
\begin{align}
&\mathbb E\int_{I}|z|^2(x,\tau){\rm
d}x\leq
e^{C(\tilde\tau-\tau)}\mathbb E\int_{I}|z|^2(x,\tilde
\tau){\rm d}x,\quad 0\leq\tau<\tilde\tau\leq T.
\end{align}
Letting $\tau=0$ and integrating over
$\left[\frac{T}{3},\frac{2T}{3}\right]$ with respect to
$\tilde\tau$, we find that
\begin{align}\label{4.28}
&\frac{T}{3}\mathbb E\int_{I}|z|^2(x,0){\rm
d}x\leq
C\mathbb E\int_{\frac{T}{3}}^{\frac{2T}{3}}\int_{I}|z|^2{\rm
d}x{\rm d}t.
\end{align}

On the other hand, by Theorem 4.2 we obtain
\begin{align}
&\mathbb E\int_{Q_T}s^3\xi^3|z|^2e^{2s\varphi}{\rm d}x{\rm d}t+\mathbb E\int_{Q_T}s\xi(x+\varepsilon)^{\alpha}|z_x|^2e^{2s\varphi}{\rm d}x{\rm d}t\nonumber\\
\leq& C\mathbb E\int_{Q_T}s^{2}\xi^{2}|Z|^2e^{2s\varphi}{\rm d}x{\rm
d}t+C\mathbb E\int_{\omega_T}s^3\xi^3|z|^2e^{2s\varphi}{\rm d}x{\rm
d}t
\end{align}
for all $\lambda\geq\lambda_3$, $s\geq s_3$. We fix
$\lambda=\lambda_3$ and $s=s_3$. By
\begin{align*}\xi^{3}e^{2s\varphi}\geq
\left(\frac{4}{T^2}\right)^{6}{\rm
exp}\left(2s_3\left(\frac{9}{2T^2}\right)^2
e^{-2\lambda_3M}\right),\quad t\in
\left[\frac{T}{3},\frac{2T}{3}\right],
\end{align*}
we further have
\begin{align}\label{4.30}
&\mathbb E\int_{\frac{T}{3}}^{\frac{2T}{3}}\int_{I}|z|^2{\rm
d}x{\rm d}t\nonumber\\
\leq &C(\lambda_3,s_3)\mathbb E\int_{Q_T}\xi^{2}|Z|^2e^{2s\varphi}{\rm d}x{\rm
d}t+
C(\lambda_3,s_3)\mathbb E\int_{\omega_T}\xi^{3}|z|^2e^{2s\varphi}{\rm d}x{\rm d}t.
\end{align}
Since $\max_{(x,t)\in Q_T}\xi^3(t)e^{2s\varphi(x,t)}<\infty$, we deduce from (\ref{4.30}) that
\begin{align}\label{4.31}
\mathbb E\int_{\frac{T}{3}}^{\frac{2T}{3}}\int_{I}|z|^2{\rm
d}x{\rm d}t\leq&
C\mathbb E\int_{Q_T}|Z|^2{\rm d}x{\rm
d}t+
C\mathbb E\int_{\omega_T}|z|^2{\rm d}x{\rm d}t.
\end{align}
Finally, we obtain the desired estimate (\ref{4.23}) from (\ref{4.28})
and (\ref{4.31}) and then complete the proof of Lemma 4.4.
\hfill$\Box$

\vspace{2mm}

Now we prove Theorem 4.1.

\vspace{2mm}

 {\noindent\bf Proof of Theorem 4.1.}\ The proof is based on a classical dual argument and an approximate method.  We introduce a linear subspace of $L^2_\mathcal F(0,T;L^2(\omega))\times L^2_\mathcal F(0,T;L^2(I))$:
 \begin{align*}
 X=\left\{(z|_\omega, Z)\ |\ (z,Z) \ {\rm solves\ the\ system\ (\ref{4.3})\ with\ some} \ z^T\in L^2(\Omega,\mathcal F_T, \mathbb P; L^2(I))\right\}
 \end{align*}
 endowed with the norm
 \begin{align*}
 \|(z|_\omega, Z)\|^2_{X}=\int_{\omega_T}|z|^2{\rm d}x{\rm d}t+\int_{Q_T} |Z|^2{\rm
d}x{\rm d}t.
 \end{align*}
  We further define a linear functional on $X$ as follows:
  \begin{align*}
  \mathcal L(z|_\omega, Z) =-\mathbb E\int_{I} y^\varepsilon(x,0)z(x,0){\rm d}x.
  \end{align*}
  By Lemma 4.4, we see that for any $\varepsilon\in (0,\nu)$, there exists constant $C$ independent of $\varepsilon$ such that
 \begin{align}
 |\mathcal L(z|_\omega, Z)|\leq &\Big(\mathbb E\int_{I}|y^\varepsilon(x,0)|^2{\rm d}x\big)^{\frac{1}{2}}\Big(\mathbb E\int_{I}|z(x,0)|^2{\rm d}x\Big)^{\frac{1}{2}}\nonumber\\
 \leq & C \Big(\mathbb E\int_{I}|y^\varepsilon(x,0)|^2{\rm d}x\Big)^{\frac{1}{2}}\|(z|_\omega, Z)\|_{X},
 \end{align}
 which means that $\mathcal L$ is a bounded linear functional on $X$.  We can extend $\mathcal L$ to be a bounded linear functional on $L^2_\mathcal F(0,T;L^2(\omega))\times L^2_\mathcal F(0,T;L^2(I))$ and use the same notation for this extension. Now by Riesz representation,  we know that for any $\varepsilon\in (0,\nu)$, there exists a unique pair of controls $(g^\varepsilon,G^\varepsilon)\in L^2_\mathcal F(0,T;L^2(\omega))\times L^2_\mathcal F(0,T;L^2(I))$ such that
  \begin{align}
  -\mathbb E\int_{I} y^\varepsilon_0(x)z(x,0){\rm d}x=\mathbb E\int_{\omega_T} g^\varepsilon z{\rm d}x{\rm d}t+\mathbb E\int_{Q_T} G^\varepsilon Z{\rm d}x{\rm d}t,
  \end{align}
  and
  \begin{align}
 \|(g^\varepsilon,G^\varepsilon)\|_{X}\leq C \Big(\mathbb E\int_{I}|y^\varepsilon_0(x)|^2{\rm d}x\Big)^{\frac{1}{2}}.
 \end{align}
   By the duality between $z$ and $y^\varepsilon$
   \begin{align}
  \mathbb E\int_{I} y^\varepsilon(x,T)z(x,T){\rm d}x -\mathbb E\int_{I} y^\varepsilon_0(x)z(x,0){\rm d}x=\mathbb E\int_{Q_T}(g^\varepsilon \mathbf 1_{\omega} z+ G^\varepsilon Z){\rm d}x{\rm d}t,
   \end{align}
  we see that for any $\varepsilon\in (0,\nu)$, there exists a pair of controls $(g^\varepsilon,G^\varepsilon)\in L^2_\mathcal F(0,T;L^2(\omega))\times L^2_\mathcal F(0,T;L^2(I))$ such that $y^\varepsilon(x,T)=0$, $\mathbb P$-$a.s$. Since the equation (\ref{4.1}) is linear,  we could further obtain that $\{(g^\varepsilon,G^\varepsilon)\}$ is Chauchy sequence such that
  \begin{align}\label{4.34}
 \|(g^{\varepsilon_1}-g^{\varepsilon_2},G^{\varepsilon_1}-G^{\varepsilon_2})\|_{X}\leq C \Big(\mathbb E\int_{I}\big|\big(y^{\varepsilon_1}_0-y^{\varepsilon_2}_0\big)(x)\big|^2{\rm d}x\Big)^{\frac{1}{2}}
 \end{align}
   for any $\varepsilon_1,\varepsilon_2\in (0,\nu)$. Notice that the constant $C$ in (\ref{4.34}) is independent of $\varepsilon$. Therefore together with $y^\varepsilon_0\rightarrow y_0$ in $L^2(\Omega,\mathcal F_0,\mathbb P; L^2(I))$, letting $\varepsilon\rightarrow 0$, we obtain a  control $(g,G)\in L^2_\mathcal F(0,T;L^2(\omega))\times L^2_\mathcal F(0,T;L^2(I))$ that drives the corresponding solution
$y$ to zero at time $T$. This completes the proof of
Theorem 4.1.\hfill$\Box$

\section{Stability for inverse problem}
\setcounter{equation}{0}

In this section, we apply Carleman estimate (\ref{3.45}) to prove the Lipschitz stability for our inverse random source problem, i.e. the following Theorem 5.1.

\begin{thm} Let $\alpha\in \left(0,2\right)$,  $r\in L^\infty_\mathcal F(0,T; W^{1,\infty}(I))$ such that $|r(x,t)|\geq r_0>0$ for $(x,t)\in Q_T$, $\mathbb P-a.s.$, $h^{(i)}\in L^2_\mathcal F(0,T)$ for $i=1,2$. Then  there exists a positive constant $C=C(\omega,I,T,\alpha,r_0)$ such that
\begin{align}\label{5.1}
&\|h^{(1)}-h^{(2)}\|_{L^2_\mathcal F(0,T)}\nonumber\\
\leq& C\left(\|y^{(1)}-y^{(2)}\|_{L^2_\mathcal F(0,T;L^2(\omega))}+\|\big(y^{(1)}-y^{(2)}\big)(\cdot, T)\|_{L^2(\Omega,\mathcal F_T,\mathbb P;L^2(I))}\right),
\end{align}
where $y^{(i)}$ is the solutions to (\ref{1.4}) corresponding to $h^{(i)}$ for $i=1, 2$, respectively.
\end{thm}

{\noindent\bf Proof.}\ Letting $\tilde y=y^{(1)}-y^{(2)}$ and $\tilde h=h^{(1)}-h^{(2)}$, we have
\begin{align}\left\{
\begin{array}{ll}
{\rm d}\tilde y-\left(x^\alpha \tilde y_x\right)_x{\rm d}t= \tilde h(t)r(x,t){\rm d}B(t),&(x,t)\in Q_T,\\
  \tilde y(1,t)=0,& t\in (0,T),\\
 {\rm and}\ \left\{\begin{array}{ll}
\tilde y(0,t)=0&{\rm for}\ \alpha\in (0,1),\\
\big(x^\alpha \tilde y_x\big)(0,t)=0 &{\rm for}\ \alpha\in [1,2),
\end{array} \right.
& t\in (0,T),\\
 \tilde y(x,0)=0,&x\in I.
 \end{array}\right.
\end{align}
Then applying (\ref{3.45}) to $\tilde y$, we obtain
\begin{align}\label{5.3}
&\mathbb E\int_{Q_T}s\lambda\Phi|\tilde h r|^2e^{2s\Phi}{\rm d}x{\rm d}t\nonumber\\
\leq & C\mathbb E\int_{Q_T} s\Phi|\tilde h r_{x}|^2e^{2s\Phi}{\rm d}x{\rm d}t+C(\lambda)\mathbb E\int_{\omega_T}s^3\Phi^3|\tilde y|^2e^{2s\Phi}{\rm d}x{\rm d}t\nonumber\\
&+C(\lambda)s^2 e^{C(\lambda)s}\|\tilde y(\cdot, T)\|_{L^2(\Omega,\mathcal F_T,\mathbb P;L^2(I))}^2.
\end{align}
By means of $|r(x,t)|\geq r_0>0$ for $(x,t)\in Q_T$, $\mathbb P-a.s.$ and choosing $\lambda$ sufficiently large to absorb the first term on the right-hand side of (\ref{5.3}), we have \begin{align}\label{5.4}
&\mathbb E\int_{Q_T}s\lambda\Phi|\tilde h|^2e^{2s\Phi}{\rm d}x{\rm d}t\nonumber\\
\leq & C(\lambda)\mathbb E\int_{\omega_T}s^3\Phi^3|\tilde y|^2e^{2s\Phi}{\rm d}x{\rm d}t+C(\lambda)s^2 e^{C(\lambda)s}\|\tilde y(\cdot, T)\|_{L^2(\Omega,\mathcal F_T,P;L^2(I))}^2.
\end{align}
Finally, using $0<\Phi e^{2s\Phi}<C(\lambda,s)$ due to the regular weight function, we deduce (\ref{5.1}) from (\ref{5.4}) and complete the proof of Theorem 5.1. \hfill$\Box$

\section{Appendix}
\renewcommand{\theequation}{A.\arabic{equation}}
\setcounter{equation}{0}
Here, we  prove Lemma 3.4 and Lemma 4.3.

\vspace{2mm}

\noindent{\bf Proof of Lemma 3.4.}\ Let $\rho_1\in C^2(\overline I)$ be a cut-function such that $0\leq \rho_1(x)\leq 1$ for $x\in  I$, $\rho_1(x)\equiv1$ for $x\in\omega^{(1)}$ and $\rho_1(x)\equiv0$ for $x\in I\setminus \omega$. By using It\^{o} formula, $({\rm d}u^\varepsilon)^2=F_1^2{\rm d}t$ and the equation of $u^\varepsilon$, we have
\begin{align}\label{A.1}
{\rm d}\left[\rho_1\xi  (u^\varepsilon)^2e^{2s\varphi}\right]=&2\rho_1\xi u^\varepsilon e^{2s\varphi}{\rm d}u^\varepsilon+\rho_1(\xi e^{2s\varphi})_t (u^\varepsilon)^2{\rm d}t+\rho_1\xi e^{2s\varphi} ({\rm d}u^\varepsilon)^2\nonumber\\
=&2\rho_1 \xi u^\varepsilon e^{2s\varphi}\left[-\left((x+\varepsilon)^\alpha u^\varepsilon_x\right)_x{\rm d}t+f_1{\rm d}t+F_1{\rm d}B(t)\right]\nonumber\\
&+\rho_1(2s\xi\varphi_t+\xi_t) (u^\varepsilon)^2e^{2s\varphi}{\rm d}t+\rho_1\xi |F_1|^2e^{2s\varphi}{\rm d}t.
\end{align}
Then integrating both side of (\ref{A.1}) in $Q_T$ and taking mathematical expectation in $\Omega$, we find
\begin{align}\label{A.2}
&2\mathbb E \int_{Q_T}\rho_1 \xi (x+\varepsilon)^\alpha |u^\varepsilon_x|^2e^{2s\varphi}{\rm d}x{\rm d}t\nonumber\\
=&2\mathbb E \int_{Q_T} \xi\left[(\rho_1e^{2s\varphi})_x(x+\varepsilon)^\alpha\right]_x |u^\varepsilon|^2{\rm d}x{\rm d}t-2\mathbb E \int_{Q_T}\rho_1\xi  u^\varepsilon f_1 e^{2s\varphi}{\rm d}x{\rm d}t\nonumber\\
&-\mathbb E \int_{Q_T}\rho_1(2s\xi\varphi_t+\xi_t) |u^\varepsilon|^2e^{2s\varphi}{\rm d}x{\rm d}t-\mathbb E \int_{Q_T}\rho_1\xi |F_1|^2e^{2s\varphi}{\rm d}x{\rm d}t\nonumber\\
\leq &C(\lambda)\mathbb E \int_{\omega_T}s^2\xi^3|u^\varepsilon|^2e^{2s\varphi}{\rm d}x{\rm d}t+C\mathbb E\int_{Q_T} s^{-2}|f_1|^2e^{2s\varphi}{\rm d}x{\rm d}t\nonumber\\
&+C\mathbb E \int_{Q_T} \xi|F_1|^2e^{2s\varphi}{\rm d}x{\rm d}t.
\end{align}
Here we have used $\left|\left[(\rho_1e^{2s\varphi})_x(x+\varepsilon)^\alpha\right]_x\right|\leq C(\lambda)s^2\xi^2 e^{2s\varphi}$ in $\omega^{(1)}$, where $C$ is independent of $\varepsilon$.
Noting that $\rho_1\equiv 1$ in $\omega^{(1)}$, we immediately deduce (\ref{3.7}) from (\ref{A.2}) and complete the proof of Lemma 3.4. \hfill $\Box$

\vspace{2mm}

Now we prove Lemma 4.3, whose proof is similar to the one in [\ref{Liu2014ECOCV}] or [\ref{Yan2018JMAA}]. Different from those papers, here we need that the estimate (\ref{4.6}) is not depending on $\varepsilon$, which is important to study our null controllability. So that we list a detailed process here.

\vspace{2mm}

\noindent{\bf Proof of Lemma 4.3.}\ As [\ref{Yan2018JMAA}], for any $\tau>0$  we set $$\varphi_\tau(x,t)=\frac{\chi(x)\psi_1(x)+(1-\chi(x))\psi_2(x)}{(t+\tau)^2(T-t+\tau)^2}$$
and \begin{align*}
\mathcal U=\Big\{(h,H)\ |\ \mathbb E\int_{\omega_T} s^{-3}\xi^{-3} |h|^2 e^{-2s\varphi}{\rm d}x{\rm d}t+\mathbb E\int_{Q_T}s^{-2}\xi^{-2} |H|^2e^{-2s\varphi}{\rm d}x{\rm d}t<\infty\Big\}.
\end{align*}
Then we consider the following constrained extremal problem
\begin{align}\label{A.3}
\mathcal J=&\frac{1}{2}\min\limits_{(h,H)\in \mathcal{U}}\bigg(\mathbb E\int_{\omega_T} s^{-3}\xi^{-3} |h|^2 e^{-2s\varphi}{\rm d}x{\rm d}t+\mathbb E\int_{Q_T}s^{-2}\xi^{-2} |H|^2e^{-2s\varphi}{\rm d}x{\rm d}t\nonumber\\
&+\mathbb E\int_{Q_T} |w|^2e^{-2s\varphi_\tau}{\rm d}x{\rm d}t+\frac{1}{\tau}\mathbb E\int_{I} |w(x,T)|^2{\rm d}x \bigg),
\end{align}
where $w$ is the solution of (\ref{4.5}) corresponding to $(h,H)$.
By the variational method in [\ref{Lions1971}], we see that for any given $\tau$, the control problem (\ref{A.3}) admits a unique optimal solution $( h_\tau, H_\tau)\in\mathcal U$ such that
\begin{align}\label{A.4}
h_\tau=-s^3\xi^3 p_\tau e^{2s\varphi},\quad H_\tau=-s^{2}\xi^{2} P_\tau e^{2s\varphi},
\end{align}
where $(p_\tau, P_\tau)$ is the solution of the following backward stochastic equation
\begin{align}\label{A.5}\left\{\begin{array}{ll}
{\rm d}p_\tau+\left((x+\varepsilon)^\alpha p_{\tau,x}\right)_x{\rm d}t=-w_\tau e^{-2s\varphi_\tau}{\rm d}t+P_\tau{\rm d}B(t),&(x,t)\in Q_T,\\
p_\tau(0,t)=p_\tau(1,t)=0,&t\in (0,T),\\
p_\tau(x,T)=\frac{1}{\tau} w_\tau(x,T), &x\in\Omega.\end{array}\right.
\end{align}
where $w_\tau \in \mathcal H$ is the solution of (\ref{4.5}) corresponding to $(h_\tau, H_\tau)$.

 Now we prove a uniform estimate for $(w_\tau, h_\tau, H_\tau)$ in $\tau$ and $\varepsilon$.  By It\^{o} formula, (\ref{4.5}), (\ref{A.4}), (\ref{A.5}) and Young's inequality, we find that
\begin{align}\label{A.6}
&\mathbb E\int_{I} p_\tau(x,T) w_\tau(x,T){\rm d}x\nonumber\\
=&\mathbb E\int_{I} p_\tau(x,0) w_\tau(x,0){\rm d}x+\mathbb E\int_{Q_T}p_\tau\left[\left((x+\varepsilon)^\alpha w_{\tau,x}\right)_x{\rm d}t+\left(s^3\xi^3 z e^{2s\varphi}+ h_\tau{\mathbf 1}_{\omega}\right){\rm d}t\right]{\rm d}x\nonumber\\
&+\mathbb E\int_{Q_T}w_\tau\left[-\left((x+\varepsilon)^\alpha p_{\tau,x}\right)_x{\rm d}t-w_\tau e^{-2s\varphi_\tau}{\rm d}t\right]{\rm d}x+\mathbb E\int_{Q_T}P_\tau H_\tau{\rm d}x{\rm d}t\nonumber\\
\leq &\epsilon \mathbb E \int_{Q_T}s^3\xi^3 |p_\tau|^2e^{2s\varphi}{\rm d}x{\rm d}t+C(\epsilon)\mathbb E \int_{Q_T}s^3\xi^3 |z|^2e^{2s\varphi}{\rm d}x{\rm d}t-\mathbb E\int_{\omega_T} s^3\xi^3 |p_\tau|^2e^{2s\varphi}{\rm d}x{\rm d}t\nonumber\\
&-\mathbb E\int_{Q_T} |w_\tau|^2e^{-2s\varphi_\tau}{\rm d}x{\rm d}t-\mathbb E \int_{Q_T}s^{2}\xi^{2} |P_\tau|^2e^{2s\varphi}{\rm d}x{\rm d}t.
\end{align}
On the other hand, since $\alpha\in \big(0,\frac{1}{2}\big)$, we can choose $\beta=\frac{4-2\alpha}{3}\in (1,2-\alpha)$ such that $2\alpha+3\beta-4=0$. Then applying Theorem 3.1 to $p_\tau$ yields that
\begin{align}\label{A.7}
\mathbb E \int_{Q_T}s^3\xi^3 |p_\tau|^2e^{2s\varphi}{\rm d}x{\rm d}t\leq &C\mathbb E\int_{Q_T} |w_\tau|^2e^{-4s\varphi_\tau+2s\varphi}{\rm d}x{\rm
d}t+C\mathbb E\int_{Q_T}
s^{2}\xi^{2}|P_\tau|^2e^{2s\varphi}{\rm d}x{\rm
d}t\nonumber\\
&+C\mathbb E\int_{\omega_T} s^3\xi^3|p_\tau|^2e^{2s\varphi}{\rm
d}x{\rm d}t.
\end{align}
Together with $\varphi_\tau\geq \varphi$, we deduce from (\ref{A.6}) and (\ref{A.7}) that
\begin{align}\label{A.8}
&\frac{1}{\tau}\mathbb E\int_{I} |w_\tau|^2(x,T){\rm d}x+\mathbb E\int_{Q_T} |w_\tau|^2e^{-2s\varphi_\tau}{\rm d}x{\rm d}t\nonumber\\
&+ E\int_{\omega_T} s^{-3}\xi^{-3} |h_\tau|^2e^{-2s\varphi}{\rm d}x{\rm d}t+\mathbb E \int_{Q_T}s^{-2}\xi^{-2} |H_\tau|^2e^{-2s\varphi}{\rm d}x{\rm d}t\nonumber\\
\leq &C \mathbb E \int_{Q_T}s^3\xi^3 |z|^2e^{2s\varphi}{\rm d}x{\rm d}t,
\end{align}
if we choose $\epsilon$ sufficiently small. Notice that
\begin{align}\label{A.9}
{\rm d}\left(s^{-2}\xi^{-2}|w_\tau|^2e^{-2s\varphi_\tau}\right)=&s^{-2}\left(\xi^{-2}e^{-2s\varphi_\tau}\right)_t|w_\tau|^2{\rm d}t\nonumber\\
&+2s^{-2}\xi^{-2}w_\tau e^{-2s\varphi_\tau}{\rm d}w_\tau+s^{-2}\xi^{-2} e^{-2s\varphi_\tau}({\rm d}w_\tau)^2.
\end{align}
Integrating both side of (\ref{A.9}) in $Q_T$, taking mathematical expectation and using the equation of $w_\tau$, we than obtain
\begin{align}\label{A.10}
&2\mathbb E\int_{Q_T}s^{-2}\xi^{-2}(x+\varepsilon)^\alpha|w_{\tau,x}|^2e^{-2s\varphi_\tau}{\rm d}x{\rm d}t\nonumber\\
=&-\mathbb E\int_{I}\left[s^{-2}\xi^{-2}|w_\tau|^2e^{-2s\varphi_\tau}\right]_{t=0}^{t=T}{\rm d}x +\mathbb E\int_{Q_T}s^{-2}(\xi^{-2}e^{-2s\varphi_\tau})_t |w_\tau|^2{\rm d}x{\rm d}t\nonumber\\
&-2\mathbb E\int_{Q_T}s^{-2}\xi^{-2}(x+\varepsilon)^\alpha(e^{-2s\varphi_\tau})_x w_\tau w_{\tau,x}{\rm d}x{\rm d}t+2\mathbb E\int_{Q_T}s\xi w_{\tau}z e^{-2s\varphi_\tau+2s\varphi}{\rm d}x{\rm d}t\nonumber\\
&+2\mathbb E\int_{\omega_T}s^{-2}\xi^{-2} w_{\tau} h_\tau e^{-2s\varphi_\tau}{\rm d}x{\rm d}t+\mathbb E\int_{Q_T}s^{-2}\xi^{-2} |H_\tau|^2e^{-2s\varphi_\tau}{\rm d}x{\rm d}t\nonumber\\
\leq & C\mathbb E\int_{Q_T}|w_\tau|^2 e^{-2s\varphi_\tau}{\rm d}x{\rm d}t+\mathbb E\int_{Q_T}s^{-2}\xi^{-2}(x+\varepsilon)^\alpha|w_{\tau,x}|^2e^{-2s\varphi_\tau}{\rm d}x{\rm d}t\nonumber\\
&+C\mathbb E \int_{Q_T}s^{2}\xi^{2} |z|^2e^{-2s\varphi_\tau+4s\varphi}{\rm d}x{\rm d}t+C\mathbb E\int_{\omega_T} s^{-3}\xi^{-3} |h_\tau|^2e^{-2s\varphi}{\rm d}x{\rm d}t\nonumber\\
&+C\mathbb E \int_{Q_T}s^{-2}\xi^{-2} |H_\tau|^2e^{-2s\varphi}{\rm d}x{\rm d}t.
\end{align}
Therefore, it follows from (\ref{A.8}) and (\ref{A.10}) that
\begin{align}
&\frac{1}{\tau}\mathbb E\int_{I} |w_\tau|^2(x,T){\rm d}x+\mathbb E\int_{Q_T}|w_\tau|^2e^{-2s\varphi_\tau}{\rm d}x{\rm d}t\nonumber\\
&+\mathbb E\int_{Q_T}s^{-2}\xi^{-2}(x+\varepsilon)^\alpha|w_{\tau,x}|^2e^{-2s\varphi_\tau}{\rm d}x{\rm d}t+\mathbb E\int_{\omega_T} s^{-3}\xi^{-3} |h_\tau|^2e^{-2s\varphi}{\rm d}x{\rm d}t\nonumber\\
&+\mathbb E \int_{Q_T}s^{-2}\xi^{-2} |H_\tau|^2e^{-2s\varphi}{\rm d}x{\rm d}t\nonumber\\
\leq &C \mathbb E \int_{Q_T}s^3\xi^3 |z|^2e^{2s\varphi}{\rm d}x{\rm d}t,
\end{align}
where $C$ is independent of $\varepsilon$ and $\tau$. Then there exist $w\in \mathcal H^1$ and $(h,H)\in \mathcal U$ such that
  \begin{align*}
  (w_\tau, h_\tau, H_\tau)\rightharpoonup (w,h,H)\quad {\rm in}\ \mathcal H^1\times \mathcal U.
  \end{align*}
  As [\ref{Liu2014ECOCV}], by letting $\tau\rightarrow 0$, we can obtain (\ref{4.6}) and $w(x,T)=0$ in $I$, $P-a.s.$  This completes the proof of this lemma.\hfill$\Box$

\vskip 1.5cm {\bf Acknowledgement.}\ This work is supported by NSFC
(No.) \vskip 1cm

\newcounter{cankao}
\begin{list}
{[\arabic{cankao}]}{\usecounter{cankao}\itemsep=0cm} \centerline{\bf
References} \vspace*{0.5cm} \small

\item\label{Alabau2006} F. Alabau-Boussouira, P. Cannarsa, G. Fragnelli, Carleman estimates for degenerate parabolic operators with applications to null controllability, J. Evol.    Equ. 6 (2006) 161-204.

\item\label{Benhassi2013} E.M. Ait Benhassi, F. Ammar Khodja, A. Hajjaj, L. Maniar, Carleman
Estimates and null controllability of coupled degenerate systems,
Evol. Equ. Control The. 2 (2013) 441-459.

\item\label{Bao2013}  G. Bao, X. Xu, An inverse random source problem in qualifying the elastic modulus of nanomaterials, Inverse Probl. 29 (2013) 015006.

\item\label{Bao2010} G. Bao, S.-N. Chow, P. Li, H. Zhou, Numerical solution of an inverse medium scattering problem with a stochastic source, Inverse Probl. 26 (2010) 074014.

\item\label{Barbu2003} V. Barbu, A. R\u{a}scanu, G. Tessitore, Carleman estimate and controllability of linear stochastic heat equations, Appl. Math. Optim. 47 (2003) 97-120.

\item\label{Bellassoued2015}M. Bellassoued, J. Le Rousseau, Carleman estimates for elliptic operators with complex coefficients. Part I: boundary value problems. Journal de Math\'{e}matiques Pures et Appliqu\'{e}es 104 (2015) 657-728.

\item\label{Cannarsa2009-1} P. Cannarsa, L. De Teresa, Controllability of 1-D coupled degenerate parabolic equations, Electron. J.  Differential Equations  2009 (2009) 1-21.

\item\label{Cannarsa2005}P. Cannarsa, P. Martinez, J. Vancostenoble, Null
controllability of degenerate heat equations,  Adv. Differential
Equations 10 (2005) 153-190.

\item\label{Cannarsa2008} P. Cannarsa, P. Martinez, J. Vancostenoble, Carleman estimates for a
class of degenerate parabolic operators, SIAM J. Control Optim.  47 (2008) 1-19.

\item\label{Du2019} R. Du, Null controllability for a class of degenerate parabolic equations with the gradient terms, J. Evol. Equ. 19 (2019) 585-613.

\item\label{Du2018} R. Du, F. Xu, Null Controllability of a Coupled Degenerate System
with the First Order Terms, J Dyn. Control Syst. 24 (2018) 83-92.

\item\label{Fadili2017} M. Fadili, L. Mania, Null controllability of
$n-$coupled degenerate parabolic systems with $m-$ controls, J. Evol.    Equ. 17 (2017) 1311-1340.

\item\label{Fu2017} X. Fu, X. Liu A weighted identity for stochastic partial differential operators and its applications, J. Differential Equations 262 (2017) 3551-3582.

\item\label{Gao2014} P. Gao, Carleman estimate and unique continuation property for the linear stochastic Korteweg-de Vries equation, Bull. Aust. Math. Soc. 90 (2014) 283-294.

\item\label{Gao2018} P. Gao, Global Carleman estimates for the linear stochastic Kuramoto-Sivashinsky equations and their applications, J. Math. Anal. Appl. 464 (2018) 725-748.

\item\label{Hu1991} Y. Hu, S. Peng, Adapted solution of a backward semilinear stochastic evolution equations, Stoch. Anal. Appl. 9 (1991) 445-459.

\item \label{Imanuvilov2001} O.Y. Imanuvilov and M. Yamamoto, Carleman estimate for a parabolic equation in a Sobolev space of negative order and its applications,
in: Control of Nonlinear Distributed Parameter Systems, in: Lecture
Notes in Pure and Appl. Math., vol.218, Dekker, New York, 2001,
pp.113-137.

\item \label{Iman3} O.Y. Imanuvilov and M. Yamamoto, Carleman estimates for the non-stationary Lem\'{e}
system and application to an inverse problem, ESAIM Control
Optim. Calc. Var. 11 (2005) 1-56.

\item \label{Isakov2} V. Isakov, N. Kim, Carleman estimates with second large parameter and applications to elasticity with residual stress. Applicationes Mathematicae 35 (2008) 447-465.

\item\label{Klib3} M.V. Klibanov, Carleman estimates for global uniqueness, stability and numerical methods for coefficient inverse problems, J. Inverse Ill-Posed Prob.
21 (2013) 477-560.

\item \label{Klib2}M.V. Klibanov, A. Timonov, Carleman Estimates for Coefficient
Inverse Problems and Numerical Applications, VSP, Utrecht, 2004.

\item\label{Krylov1994} N.V. Krylov, A $W^n_2$-theory of the Dirichlet problem for SPDEs in general smooth domains. Probab. Theory Relat. Fields 98 (1994) 389-421.

\item\label{Lions1971} J.L. Lions, Optimal Control of Systems Governed by Partial Differential Equations, vol.170, Springer-Verlag, New York-Berlin, 1971.

\item\label{Liu2018MCRF} L. Liu and X. Liu, Controllability and observability of some coupled stochastic parabolic systems, Math. Control. Relat. F. 8 (2018) 829-854.

\item\label{Liu2014ECOCV}X. Liu, Global Carleman estimate for stochastic parabolic equations and its application, ESAIM: Control Optim. Calc. Var. 20 (2014) 823-839.

\item \label {LiuSIAM} X. Liu,  Y. Yu, Carleman Estimates of Some Stochastic Degenerate Parabolic Equations and Application, SIAM J. Control Optim. 57(2019) 3527-3552.

\item\label{Lv2012} Q. L\"{u}, Carleman estimate for stochastic parabolic equations and inverse stochastic parabolic problems, Inverse Probl. 28 (2012) 045008.

\item\label{Lv2013} Q. L\"{u}, Exact controllability for stochastic Schr\"{o}dinger equations, J. Differential Equations 255 (2013) 2484-2504.

\item\label{Lv2015}Q. L\"{u} and X. Zhang,  Global uniqueness for an inverse stochastic hyperbolic problem with three unknowns, Commun. Pur. Appl. Math. 68 (2015) 948-963.

\item\label{Opic1990} B. Opic and A. Kufner, Hardy-type inequalities, Halsted Press, 1990.

\item\label{Rousseau} J.L. Rousseau and G. Lebeau, On Carleman estimates for elliptic and parabolic operators.
Applications to unique continuation and control of parabolic equations, ESAIM Control
Optim. Calc. Var. 18 (2012) 712-747.

\item\label{TangSIAM2009} S. Tang and X. Zhang, Null controllability for forward and backward stochastic parabolic equations, SIAM J. Control Optim. 48 (2009) 2191-2216.

\item\label{Wang2014} C. Wang, R. Du, Carleman estimates and null controllability for a
class of degenerate parabolic equations with convection terms, SIAM J. Control Optim. 52 (2014) 1457-1480.

\item\label{Wu2019} B. Wu, Y. Gao,  Z. Wang and Q. Chen, Unique continuation for a reaction-diffusion system with cross diffusion, J. Inverse Ill-Posed Probl. 27 (2019) 511-525.

\item\label{Wu2018} B. Wu, L. Yan, Y. Gao and Q. Chen, Carleman estimate for a linearized
bidomain model in electrocardiology and its applications, Nonlinear Differ. Equ. Appl. 25 (2018) 4 (20pp).

\item\label{Yamamoto1}M. Yamamoto, Carleman estimates for parabolic equations and
applications,  Inverse Probl. 25 (2009) 123013(75pp).

\item\label{Yan2015IP} G. Yuan,  Determination of two kinds of sources simultaneously for a stochastic wave equation,  Inverse Probl. 31 (2015) 085003.

\item\label{Yuan2017JMAA} G. Yuan, Determination of two unknowns simultaneously for stochastic Euler-Bernoulli beam equations, J. Math. Anal. Appl. 450 (2017) 137-151.

\item\label{Yan2018JMAA} Y. Yan, Carleman estimates for stochastic parabolic equations with Neumann boundary conditions and applications, J. Math. Anal. Appl. 457 (2018) 248-272.

\item\label{ZhangSIAM2008} X. Zhang, Carleman and observability estimates for stochastic wave equations, SIAM J. Math. Anal. 40 (2008) 851-868.

\end{list}

\end{document}